\documentclass[11pt, a4paper]{article}%
\usepackage{amsmath}
\usepackage{amsfonts,xcolor}
\usepackage{amssymb}
\usepackage{graphicx,sw20jart}
\usepackage{caption,subcaption,eurosym}
\usepackage{graphicx}
\usepackage{comment}
\usepackage{amsfonts}%
\providecommand{\U}[1]{\protect\rule{.1in}{.1in}}
\newtheorem{theorem}{Theorem}[section]

\newtheorem{definition}{Definition}[section]

\newtheorem{lemma}[theorem]{Lemma}
\newtheorem{proposition}[theorem]{Proposition}
\newtheorem{remark}{Remark}[section]

\begin{document}

\title{Optimal Abatement Schedules for Excess Carbon Emissions\\ Towards a Net-Zero Target}%Optimal Carbon Emissions Under a Ratcheting-Down Constraint}
\author{Hansj\"{o}rg Albrecher\thanks{Department of Actuarial Science, Faculty of
		Business and Economics and Swiss
		Finance Institute and Expertise Center for Climate Extremes, University of Lausanne, CH-1015 Lausanne. Email: hansjoerg.albrecher@unil.ch} \and  Nora Muler\thanks{Departamento de Matemática y Estadística, Universidad
		Torcuato Di Tella. Av. Figueroa Alcorta 7350 (C1428BIJ) Ciudad de Buenos
		Aires, Argentina. Email: nmuler@utdt.edu}}
\date{}
\maketitle

\abstract{\begin{quote} 
		Achieving net-zero carbon emissions requires a transformation of energy systems, industrial processes, and consumption patterns. In particular, a transition towards that goal involves a gradual reduction of excess carbon emissions that are not essential for the well-functioning of society. In this paper we study this problem from a stochastic control perspective to identify the optimal gradual reduction of the emission rate, when an allocated excess carbon budget is used up over time. Assuming that updates of the available carbon budget follow a diffusion process, we identify the emission strategy that maximizes {\color{black}the expected discounted future profit from these excess emissions} under the constraint of a non-increasing emission rate, with an additional term rewarding the amount of time for which the {\color{black}excess carbon budget} is not yet depleted. We establish a link of this topic to optimal dividend problems in insurance
		risk theory under ratcheting constraints and show that the value function is the unique viscosity solution of the associated Hamilton-Jacobi-Bellman equation. We provide numerical illustrations of the resulting optimal abatement schedule of emissions and a quantitative evaluation of the effect of the non-increasing rate constraint on the value function.
\end{quote}}

\section{Introduction}
Motivated by the Paris Agreement adopted within the United Nations Framework Convention on Climate Change (UNFCCC), whose objective is to curb global anthropogenic greenhouse gas (GHG) emissions (see, for example, Popovski \cite{Popo}), many governments have recently announced commitments to reach net-zero carbon emissions by specified target dates.\footnote{Like many other authors, for simplicity we refer to GHG emissions as carbon emissions in this paper, as carbon dioxide and methane correspond to more than 90\% of the GHG emissions.} Such a goal can only be achieved by substantially reducing avoidable emissions and compensating for those that are unavoidable. Naturally, these reductions are difficult to realize in light of established consumption patterns and the significant inertia associated with behavioral change. This applies equally to individuals, firms, and society as a whole, and the political dimension of this question is at this point mainly on plans and rules for companies, and possibly the exertion of implicit and explicit pressure from governments through respective directives and laws.\\
This topic can be examined from multiple perspectives, see for instance Borissov and Bretschger \cite{Borissov} for an economic viewpoint on fair  contributions across countries with heterogeneous wealth and pollution intensity level. Once a carbon emission target is set on a country level, its implementation as a tradeoff between emission trading and actual emission abatement is a non-trivial task, see for instance A{\"\i}d and Biagini \cite{Aid2,Aid} and Biagini \cite{Biagini} for the study of this as a stochastic Stackelberg game between firms and the regulator, cf.\ also Wijk \cite{Wijk}. For an analysis and tracking of the transition path of an individual company towards a net-zero target, see for instance Chekriy et al.\ \cite{Kiesel} and Saleh et al.\ \cite{tankov2}. Huang et al.\ \cite{Huang} examine a stochastic control problem for carbon emission reduction and the purchase amount of carbon allowances as a bivariate control problem, see also Chen  et al.\ \cite{Cheng}. 

For a profit-maximizing company the tradeoff between paying carbon taxes and investing into technologies to reduce carbon emissions can lead to an interesting stochastic control problem, see Colaneri et al.\ \cite{Col}. Bourgey et al.\ \cite{Bourgey} study another dynamic control problem of maximizing profit (which is increasing with the intensity of carbon emission) when at the same time facing penalties as a function of the discrepancy between the actual emission rate and a target emission rate that reduces according to a given socio-economic pathway (SSP).  \\

In many situations, it makes sense to assume that for reaching a net-zero target over time, one decides to compensate inevitable carbon emissions through the purchase of carbon allowances or other compensation mechanisms (like investing into carbon sequestration etc.), and then is left with a budget for (a priori) avoidable 
excess emissions that can be used up until the time at which the net-zero goal should be achieved (or is politically enforced). These excess emissions could be linked to profit when considering a company, or may serve some personal utility if one has an individual in mind. The question is then the schedule according to which this excess emission budget shall optimally be consumed, and a profound understanding of such patterns could be helpful to develop appropriate incentives for successfully reaching net-zero goals. {\color{black}The excess carbon} budget may itself be subject to uncertainty over time (increases to due to technology advances in carbon capture facilities, decreases due to stronger political pressure, {\color{black}changes in regulatory standards,} etc., or simply noise), and it may be useful to model the available excess carbon emission budget as a stochastic process. {\color{black}For simplicity, we will model this excess budget in this paper by a diffusion process.} Albrecher and Zhu \cite{AZ25} recently studied such a problem with techniques from stochastic control theory. Concretely, they looked at the problem of when and at which rates to optimally use up an available excess carbon budget, if one assigns value to the carbon emissions according to a linear utility function and expresses a preference for earlier emission by using a constant discount rate. Once the budget is depleted, there will be no future excess emissions. The availability of some carbon budget at future times is rewarded by a constant term $\Lambda$ that is also subject to the same discounting and is added to the value function whenever the process is not yet depleted. Hence, $\Lambda$ can be considered to represent a certain {\color{black}monetary} sustainability component in the optimization, or also a {\color{black}quantitatvely expressed} desire to leave some excess carbon emission for later (or even for future generations); see e.g. Korn \cite{korn}, Korn and Nurkanovic \cite{korn2} for other proposals to incorporate sustainability aspects in profitability considerations. In Albrecher and Zhu \cite{AZ25}, %the carbon budget was assumed to follow a diffusion process. For that situation, 
the optimal emission strategy was identified as a barrier strategy, with maximally allowed emission rate as soon as the available carbon budget exceeds the barrier, and no emissions below that barrier. The emphasis was then on effects of present-bias (linked to subjective discount rates) on respective emission decisions, and the corresponding efficiency of carbon taxation towards the net-zero target.  

In the present paper we would like to take a different angle on identifying optimal excess carbon emission schedules. Even if a barrier strategy is optimal for maximizing the expected discounted excess emissions until depletion with a $\Lambda$-reward on keeping the budget positive, its implementation results in a lot of variability in consumption patterns, as there is no emission below the barrier and maximally allowed emission above it. It may be easier (both psychologically and practically) to implement incentives or requirements that foresee a gradual reduction of the excess carbon emission of the entity (individuals, companies or the society at large) until the excess emission budget is used up. It is therefore of interest to consider the optimal emission problem with the constraint that emission rates can only decrease, and the challenge is then to find the strategy that optimizes emissions according to the above objective under this abatement constraint. In particular, it is desirable to quantify the efficiency loss which this constraint entails. From a methodological perspective, there is a certain degree of similarity of the present stochastic control problem with identifying optimal dividend payout strategies from an insurance portfolio in classical risk theory, when the goal is to maximize expected discounted dividends until ruin (see e.g.\ Albrecher and Thonhauser \cite{AT09}). For this latter problem, in Albrecher et al.\ \cite{AAM,AAM1} a ratcheting constraint was considered, where dividend rates can never be decreased. The setup of the present paper refers to the situation where such dividend rates (excess emission rates in the present context) can, in contrast, never be increased (which we also occasionally will refer to as 'down-ratcheting' in the sequel). Certain parts of the proofs of our results correspondingly benefit from similarities to proofs that were developed in Albrecher et al.\ \cite{AAM1}. 

%
%This work presents a stochastic control model for optimal carbon-dioxide
%emission management under environmental constraints. We consider the $CO_{2}$
%surplus of a country (or company) as a stochastic process modeled by a
%Brownian motion with positive drift and volatility, representing the random
%evolution of emission reserves. The country chooses its $CO_{2}$ emission
%rates from a bounded set of admissible controls, where only ratcheting-down
%strategies are permitted, i.e. emission rates can be reduced or maintained
%over time but never increased.

We define a performance criterion that accounts for expected  {\color{black}profit due to} 
cumulative discounted excess emissions, together with a constant reward $\Lambda$ for
safeguarding unused carbon emission capacity, until the time of depletion of that excess carbon budget (the time when the
controlled surplus first becomes negative). The objective is to determine the
optimal excess emission strategy that maximizes this function under the down-ratcheting constraint, which we achieve by identifying the optimal strategy to be of threshold type for a discretized version of the problem and then showing uniform convergence of the discrete problem to the continuous one. This results into an optimal excess emission abatement schedule (or \textit{emission abatement curve}), according to which emissions are permanently reduced to a lower level whenever new record lows of the still available carbon budget have been reached, until the excess emissions are reduced to zero. We also illustrate the approach for a few concrete numerical examples with positive, zero and negative drift of the carbon excess emission surplus process, and compare the optimal emission strategy to the one without the abatement constraint as well as to the situation where one simply applies a linear reduction of the excess emission rate over time. The latter helps to see the degree of performance increase that is possible through the application of the optimal excess emission abatement schedule.

{\color{black}To summarize, previous contributions in the stochastic control literature mentioned above have mainly focused on optimal compensation mechanisms through the purchase of carbon certificates and carbon sequestration, together with their inherent uncertainties regarding cost and availability. They have generally not addressed changes in emission habits or profit objectives themselves. In this paper we aim to shed some light on the optimal timing of effective reductions in emission habits when such reductions are inevitable in the long run, but remain, at least to some extent, voluntary at the present moment. For this purpose, we distinguish between unavoidable emissions, for which compensation mechanisms are assumed to be in place, and avoidable excess emissions, for which a remaining budget is available. The objective function used to identify the optimal abatement strategy is still linked to the monetary profit generated by these remaining emissions (for individuals, a translation of the preference to emit into a monetary unit would be needed). However, we add a second term that counterbalances emission-related profit with a monetary reward for not having depleted the excess budget at any point in time. This reward is discounted at the same rate as the profit term, reflecting the same degree of preference for current profits and rewards over future ones. In particular, we are interested in determining how large such a reward term, $\Lambda$, must be in order to significantly shift the optimal strategy away from a purely profit-driven one and towards one involving abatement. After the theoretical analysis, Section 7 illustrates the results using realistic parameter magnitudes for a large firm.\\
}

The remainder of the paper is organized as follows. Section 2 introduces the
model and the detailed formulation of the problem. It also provides some first
basic results on properties of the value function under consideration. Section
3 derives the Hamilton-Jacobi-Bellman (HJB) equation and shows that the value function is a viscosity solution of the HJB equation, together with a verification theorem. In Section 4 we formulate the problem on a discrete set of admissible emission rates and in Section 5 we prove
that the value function of the problem for discrete sets convergences
to the one for a continuum of admissible emission rates as the mesh size of
the finite set tends to zero. The latter paves the way for establishing optimal solutions numerically in an efficient way. In Section 6 we show that for finitely many
admissible emission rates, there exists an optimal strategy for which the
change and non-change regions have only one connected component (this
corresponds to the extension of one-dimensional threshold strategies to the
two-dimensional case). We also provide an implicit equation defining the
optimal threshold function for this case. Section 7 then contains numerical
illustrations of the optimal strategy and comparisons to the unconstrained case as well as to the simpler strategy of linearly reducing emission rates over time. Section 8 concludes and identifies some future research directions of interest. Some technical proofs are delegated to an appendix.

\section{Model and basic results \label{Model and basic results}}

Assume that the (excess) carbon emission budget of an entity (a country, a company or even an individual person) available at time $t$ is modeled by a Brownian motion with drift:%
\begin{equation}
X_{t}=x+\mu t+\sigma W_{t},\label{BM}%
\end{equation}
where $W_{t}$ is a standard Brownian motion, and $\sigma>0,\, \mu\in{\mathbb R}$ are given
constants.\footnote{ As described above, $X_t$ refers to \textit{excess} emissions that are a priori avoidable, so in the sequel the term 'emissions' will always refer to these 'excess emissions'. The unit of $X_t$ could for e.g.\ be {\color{black}MtCO$_2$}.} The entity uses this budget $X_t$ to emit carbon at rates
chosen from a set $S\subset\lbrack0,\overline{c}]$ , where $\overline{c}\geq0$
is the maximum allowable emission rate. {\color{black}The specification \eqref{BM} is a tractable reduced-form diffusion
model for the uncertain evolution of the remaining excess emission budget. The drift $\mu$ represents the expected trend in the budget, while the
volatility $\sigma$ captures aggregate fluctuations due to economic,
technological, regulatory, measurement, and weather-related shocks (realistic parameter values for $\mu,\sigma$ as well as for $x$ vary with the entity being considered, see Section \ref{Numerical examples}). This Browniam specification is
in line with recent stochastic control models
of carbon emission markets, see, e.g., A{\"\i}d and Biagini \cite{Aid2}.}\\

Let $(\Omega,\mathcal{F},\left(  \mathcal{F}_{t}\right)  _{t\geq
	0},\mathcal{P})$ be the complete probability space generated by the process
$(X_{t})_{t\geq0}$, and let $C_{t}$ denote the {\color{black}excess} carbon emission rate at time $t$. In this paper, we want to only consider emission patterns where the {\color{black}excess} emission rate can not be increased beyond its current level anymore. 
 Given an initial budget $X_{0}=x$
and an initial emission rate $c\in S$ at $t=0$ (which typically will be equal to $\overline{c}$), an admissible strategy is therefore a
process $C=\left(  C_{t}\right)  _{t\geq0}$ that is non-increasing,
right-continuous and adapted to the filtration $\left(  \mathcal{F}%
_{t}\right)  _{t\geq0}$ with $C_{t}\in S$ for all $t$. In other words, the
{\color{black}entity (country, company or individual person)} is only allowed to reduce or maintain its {\color{black}excess} emission rate over time, so admissible emission strategies are all of the \textit{ratcheting-down} type. Under a given
strategy $C$, the controlled carbon emission surplus process can be written as%
\begin{equation}
X_{t}^{C}=X_{t}-\int_{0}^{t}C_{s}ds. \label{XtC}%
\end{equation}
Define $\Pi_{x,c}^{S}$ as the set of all admissible ratcheting-down strategies with initial surplus $x\geq0$ and initial emission 
rate $c\in S$. \ Given $C\in\Pi_{x,c}^{S}$, the value function of this
strategy that we consider in this paper includes a reward for not having exhausted the remaining budget too early and is given by%
\begin{equation}
J_{\gamma}(x;C)=\mathbb{E}\left[  \int_{0}^{\tau}e^{-qs}({\color{black}\gamma\,C_{s}+\Lambda_{\gamma}})ds\right] .\label{Objfunction}%
\end{equation}
{\color{black}Let us describe this choice is more detail. We assume that the monetary payoff associated with emitting at rate $C_s$ is proportional to the emission rate, with proportionality constant $\gamma=c_{net}-c_{tax}>0$. Here, $c_{net}$ denotes the net profit per unit rate of emissions, for example the revenue from selling the produced good minus the physical costs of its production, expressed per unit of associated emissions (if the considered entity is an individual or society as a whole, $c_{net}$ may instead be interpreted more broadly as a monetary valuation of the benefit associated with the corresponding emission-generating activity). The term $c_{tax}$ represents a per-unit emission charge, such as a carbon tax, which may be interpreted as reflecting the social cost of carbon emissions. The quantity $\tau=\inf\left\{  t\geq0:X_{t}^{C}<0\right\}$ is the depletion time,
i.e., the first time the controlled {\color{black}remaining excess budget} becomes negative. At the same time, the term $\Lambda_{\gamma}>0$ is a fixed reward parameter that assigns a monetary value to maintaining a positive excess emission budget. It therefore makes the trade-off between current profits from emissions and preserving a remaining emission budget quantitatively comparable. Finally, $q>0$ denotes the discount rate.\\
Hence, for any
initial {\color{black}excess budget} $x\geq0$ and initial emission rate $c$, our aim is to maximize {\color{black}the optimal value function}%

\begin{equation}
V^{S}(x,c)=\sup_{C\in\Pi_{x,c}^{S}}J_{\gamma}(x;C)=\gamma \sup_{C\in\Pi_{x,c}^{S}}J_1(x;C), \label{Optimal Value Function}%
\end{equation}
with $\Lambda_1:=\Lambda=\Lambda_{\gamma}/\gamma$. That is, for every value of $\gamma$ we can w.l.o.g.\ reduce the optimization problem to $\gamma=1$ by dividing the original parameter $\Lambda_{\gamma}$ by $\gamma$, which we will do in the sequel.  

% {\color{black}\begin{remark}
% \normalfont  In the absence of uncertainty about the available budget in the future (i.e., $\sigma=0$ in \eqref{BM}), without such a $\Lambda$-term there would be no incentive to reduce emissions below the maximally allowed emission rate $\overline{c}$. With $\sigma>0$, the decision-maker may benefit from budget increases  \end{remark} }

\begin{remark}
\normalfont Note that in \eqref{Optimal Value Function} we do not explicitly account for potential costs of emission abatement. The rationale is that $X_t$ represents only  the budget of a priori avoidable excess emissions. Reducing these emissions is therefore assumed to be feasible, although potentially inconvenient, and this inconvenience is reflected in the resulting reduction of the profit rate $\gamma\,C_{s}$. One may interpret reductions in the emission rate $C_s$ as being implemented by first exploiting the least costly abatement options and only subsequently moving toward more expensive ones. Under sufficient anticipation, such abatement costs may thus be kept small, or even negligible, relative to the profit loss already captured by the lower emission rate.\\
In addition, the model does not explicitly allow for the purchase of carbon certificates to increase the budget $X_t$. This is consistent with the interpretation of $X_t$ as an excess budget for emissions that are, in principle, avoidable. Unavoidable emissions are assumed to have been accounted for separately, either by subtracting them from the initial budget $X_0$, by incorporating them into the drift $\mu$ of the process $X_t$, or by covering them through carbon certificates outside the model. So the focus of this paper is on how to best reduce excess emissions when preferences can be quantified through the above combination of $\gamma\,C_{s}$ and $\Lambda_{\gamma}$, under the intuitive requirement that $C_{s}$ is decreasing over time.
\end{remark}}

\begin{remark}
\normalfont\label{Optima sin ratcheting} {\color{black} If the carbon emission budget process \eqref{BM} were interpreted as a surplus process of an insurance portfolio and the emissions as dividend payments, then our optimal value function \eqref{Optimal Value Function} corresponds to the one of an optimal dividend problem (maximizing expected discounted bounded dividends until ruin), where a delayed time of ruin is rewarded through a monetary rate $\Lambda$ as long as the process remains positive. Such a problem was indeed studied in Thonhauser and Albrecher \cite{ThAl07} under the dynamics \eqref{BM} without the ratcheting-down constraint on the payout strategy considered in this paper. Denote the value function of this one-dimensional dividend control problem for the same parameters $\mu,\sigma$ and $\Lambda$ by $V_D(x)$. While the ratcheting-down constraint is not of immediate interest in the traditional dividend context, this connection and therefore the results of \cite{ThAl07} provide an upper bound for our present value function $V^{S}(x,c)$: }
%Concretely, the unconstrained optimal dividend problem in such 
%with a reward for a later time of ruin through the same factor $\Lambda$ in the objective function as in \eqref{Objfunction} (with $\gamma=1$), which therefore provides an upper bound for the value function in this paper.
%The unconstrained dividend problem is one-dimensional;
%let us denote its value function for the same parameters $\mu,\sigma$ and $\Lambda$ by $V_{D}(x)$. 
We have $V^{S}(x,c)\leq
V_{D}(x)$ for all $x\geq0$ and $c\in S\subset\lbrack0,\overline{c}]$. The
function $V_{D}$ is increasing, concave, twice continuously differentiable
with $V_{D}(0)=0$, $\lim_{x\rightarrow\infty}V_{D}(x)=(\overline{c}%
+\Lambda)/q$ and $V_{D}^{\prime}(x)\leq V_{D}^{\prime}(0)$ for all $x\geq0$. \hfill $\diamond$
\end{remark}

\begin{remark}
\normalfont\label{Optima con ratcheting}Our optimal stochastic control problem
is also related to the classical dividend optimization problem with a
ratcheting-up constraint in insurance surplus models (see Albrecher et al.\ \cite{AAM1} and Guan and Xu \cite{GQ}). However, in contrast to these works, we consider here a ratcheting-down constraint and incorporate the reward term
$\Lambda$.\hfill $\diamond$
\end{remark}

From the Brownian motion assumption, it is immediate that $V^{S}(0,c)=0$ for all $c\in S$, reflecting the
fact that no emissions can be sustained once the {\color{black}budget} is depleted. \\
We next establish a basic result concerning the boundedness and monotonicity
properties of the optimal value function.

\begin{proposition}
\label{Monotone Optimal Value Function}The optimal value function $V^{S}(x,c)$
is bounded above by $(\overline{c}+\Lambda)/q$, and it is non-decreasing in
both the surplus $x$ and the emission rate $c.$
\end{proposition}

\noindent\textit{Proof.} Since
\[
V^{S}(x,c)\leq V_{D}(x)\leq\frac{\overline{c}+\Lambda}{q},%
\]
we have the boundedness result.

To show monotonicity in $c$, note that if $c_{1}<c_{2}$ then $\Pi_{x,c_{1}%
}^{S}\subset\Pi_{x,c_{2}}^{S}$ for any $x\geq0$ and so $V^{S}(x,c_{1})\leq
V^{S}(x,c_{2})$.

For monotonicity in $x$, consider $0\leq x_{1}<x_{2}$ and an admissible
ratcheting-down strategy $C^{1}\in\Pi_{x_{1},c}^{S}$ for any $c\in S$, and let us
define $C^{2}\in\Pi_{x_{2},c}^{S}$ as $C_{t}^{2}=C_{t}^{1}$ until the
exhaustion time of the controlled process $X_{t}^{C^{1}}$, and then setting
$C_{t}^{2}=0$ (i.e. no emissions) afterwards. \ Clearly, $J(x;C_{1})\leq
J(x;C_{2})$ and so $V^{S}(x_{1},c)\leq V^{S}(x_{2},c)$.\hfill $\blacksquare$\newline

The following proposition provides a global Lipschitz estimate for the optimal
value function. The proof is identical to the one of Proposition 2.2 in Albrecher et al.\ \cite{AAM1}, with the obvious adaptations for the factor $\Lambda$.

\begin{proposition}
\label{Proposition Global Lipschitz zone}There exists a constant $K>0$ such
that
\[
0\leq V^{S}(x_{2},c_{1})-V^{S}(x_{1},c_{2})\leq K\left[  \left(  x_{2}%
-x_{1}\right)  +\left(  c_{2}-c_{1}\right)  \right]
\]
for all $0\leq x_{1}\leq x_{2}$ and $c_{1},c_{2}\in S$ with $c_{1}\leq c_{2}.
$
\end{proposition}

Finally, we state the Dynamic Programming Principle (DPP), its proof is
similar to the one of Lemma 1.2 in Azcue and Muler \cite{AM Libro}{\Large .}

\begin{lemma}
\label{DPP} Given any stopping time $\widetilde{\tau}$, we can write
\[
V^{S}(x,c)=\sup\limits_{C\in\Pi_{x,c}^{S}}\mathbb{E}\left[  \int_{0}%
^{\tau\wedge\widetilde{\tau}}e^{-qs}(C_{s}+\Lambda)ds+e^{-q(\tau
\wedge\widetilde{\tau})}V^{S}(X_{\tau\wedge\widetilde{\tau}}^{C},C_{\tau
\wedge\widetilde{\tau}})\right]  \text{.}%
\]

\end{lemma}

{ }

\section{Hamilton-Jacobi-Bellman equations
\label{Hamilton-Jacobi-Bellman equations}}

In this section, we introduce the HJB equation
associated with the ratcheting-down emission control problem where
the set of possible emission rates is $S:=[0,\overline{c}]\subset\lbrack
0,\infty)$ with $\overline{c}>0$. We show that the optimal value function {\color{black}$V^S,$}
defined in (\ref{Optimal Value Function}), is the unique viscosity solution of
the corresponding HJB equation with boundary condition $(\overline{c}%
+\Lambda)/q$ when $x$ goes to infinity.

Consider the strategy that emits at a constant rate $c$ until the carbon budget is exhausted. The corresponding value function $W^{c}(x)$ is the
unique solution of the second-order differential equation%
\begin{equation}
\mathcal{L}^{c}(W):=\frac{\sigma^{2}}{2}\partial_{xx}W+(\mu-c)\partial
_{x}W-qW+c+\Lambda=0 \label{Lc}%
\end{equation}
with boundary conditions $W^{c}(0)=0$ and $\lim_{x\rightarrow\infty}$
$W^{c}(x)=(c+\Lambda)/q,$ {\color{black}which is a classical problem of a first-passage-time of Brownian motion with drift, see e.g.\ Kyprianou \cite{Kyp}. Recall that} the general solutions $\mathcal{L}^{c}(W)=0$ of this differential equation are of
the form
\begin{equation}
\frac{c+\Lambda}{q}+a_{1}e^{\theta_{1}(c)x}+a_{2}e^{\theta_{2}(c)x}\text{ with
}a_{1},a_{2}\in{\mathbb{R}}, \label{Solucion General L=0}%
\end{equation}
where $\theta_{1}(c)<0<\theta_{2}(c)$ are the roots of the characteristic
equation:
\[
\frac{\sigma^{2}}{2}z^{2}+(\mu-c)z-q=0
\]
associated to the operator $\mathcal{L}^{c}$, and so%
\begin{equation}
\theta_{1}(c):=\frac{c-\mu-\sqrt{(c-\mu)^{2}+2q\sigma^{2}}}{\sigma^{2}}%
,\quad\theta_{2}(c):=\text{$\frac{c-\mu+\sqrt{(c-\mu)^{2}+2q\sigma^{2}}%
}{\sigma^{2}}$.} \label{Definicion tita1 tita2}%
\end{equation}

Since the value function must remain bounded we can discard the exponentially
growing term and the bounded solutions can be written as

\begin{equation}
\frac{c+\Lambda}{q}+ae^{\theta_{1}(c)x}\ \text{with }a\in{\mathbb{R}}.
\label{Bounded solutions}%
\end{equation}
From the boundary conditions, we then get 
\begin{equation}
W^{c}(x)=\frac{c+\Lambda}{q}\left(  1-e^{\theta_{1}(c)x}\right)  .
\label{Formula de V con c constante}%
\end{equation}
It
follows that $W^{c}(x)$ is increasing and concave.

\begin{remark}
\normalfont\label{lim l0/c}Given a set $S:=[0,\overline{c}]\subset
\lbrack0,\infty)$ $,$ we have that
\[
\frac{\overline{c}+\Lambda}{q}\geq V^{S}(x,c)\geq W^{\overline{c}%
}(x)=V^{\left\{  \overline{c}\right\}  }(x,\overline{c})=\frac{\overline
{c}+\Lambda}{q}\left(  1-e^{\theta_{1}(\overline{c})x}\right)
\]
and, consequently, $\lim_{x\rightarrow\infty}V^{S}(x,c)=(\overline{c}%
+\Lambda)/{q}$ for any $c\in S$.\hfill $\diamond$
\end{remark}

We now consider the general case where the admissible emission set is $S=[0,\overline
{c}]$ for some $\overline{c}>0.$ The HJB equation associated to
(\ref{Optimal Value Function}) is given by%

\begin{equation}
\max\{\mathcal{L}^{c}(u)(x,c),-\partial_{c}u(x,c)\}=0\text{ for }%
x\geq0\ \text{and }0\leq c\leq\overline{c}\text{,} \label{HJB equation}%
\end{equation}
where $\mathcal{L}^{c}$ is defined in (\ref{Lc}).

We say that a function $f:[0,\infty)\times\lbrack0,\overline{c})\rightarrow
{\mathbb{R}}$ is \textit{(2,1)-differentiable} if $f$ is continuously
differentiable and $\partial_{x}f(\cdot,c)$ is continuously differentiable for
each $c\in\lbrack0,\overline{c})$. To solve the HJB equation, we work in the
framework of viscosity solutions.

\begin{definition}
\label{Viscosity}

(a) A locally Lipschitz function $\overline{u}:[0,\infty)\times\lbrack
0,\overline{c}]\rightarrow{\mathbb{R}}$\ \ is a viscosity supersolution of
(\ref{HJB equation})\ at $(x,c)\in(0,\infty)\times\lbrack0,\overline{c})$\ if
any (2,1)-differentiable function $\varphi:[0,\infty)\times\lbrack
0,\overline{c}]\rightarrow{\mathbb{R}}\ $with $\varphi(x,c)=\overline{u}(x,c)$, 
and such that $\overline{u}-\varphi$ reaches the minimum at $\left(
x,c\right)  $, satisfies
\[
\max\left\{  \mathcal{L}^{c}(\varphi)(x,c),-\partial_{c}\varphi(x,c)\right\}
\leq0.\
\]
The function $\varphi$ is called a \textbf{test function for supersolution} at
$(x,c)$.

(b) A locally Lipschitz function $\underline{u}:$ $[0,\infty)\times
\lbrack0,\overline{c}]\rightarrow{\mathbb{R}}\ $ \ is a viscosity subsolution
of (\ref{HJB equation})\ at $(x,c)\in(0,\infty)\times\lbrack0,\overline{c}%
)$\ if any (2,1)-differentiable function $\psi:[0,\infty)\times
\lbrack0,\overline{c}]\rightarrow{\mathbb{R}}\ $with $\psi(x,c)=\underline
{u}(x,c)$, and such that $\underline{u}-\psi$\ reaches the maximum at $\left(
x,c\right)  $, satisfies
\[
\max\left\{  \mathcal{L}^{c}(\psi)(x,c),-\partial_{c}\psi(x,c)\right\}
\geq0\text{.}%
\]
The function $\psi$ is called a \textbf{test function for subsolution} at
$(x,c)$.

(c) A function $u:[0,\infty)\times\lbrack0,\overline{c}]\rightarrow
{\mathbb{R}}$ which is both a supersolution and subsolution at $(x,c)\in
\lbrack0,\infty)\times\lbrack0,\overline{c})$ is called a viscosity solution
of (\ref{HJB equation})\ at $(x,c)$.
\end{definition}

\begin{remark}
\normalfont\label{No importa el c}In order to simplify the notation, we define
$V(x,c):=V^{S}(x,c)$. Because of the ratcheting-down constraint on the
emission rate, we have{ }$V^{S}(x,c)=V^{[0,c]}(x,c).$\hfill $\diamond$
\end{remark}

We first prove that $V$ is a viscosity solution of the corresponding HJB
equation. The proof is in the appendix.

\begin{proposition}
\label{Proposicion Viscosidad} $V$ is a viscosity solution of
(\ref{HJB equation}) in $(0,\infty)\times\lbrack0,\overline{c}]$.
\end{proposition}

When $c=0$, the ratcheting-down constraint implies that the emissions are stopped.
Hence, $V(x,0)$ corresponds to the value function of the strategy that does
not emit, with initial surplus $x$. So, by (\ref{Formula de V con c constante}%
),
\begin{equation}
V(x,0)=V^{\{0\}}(x,0)=\frac{\Lambda}{q}\left(  1-e^{\theta_{1}(0)x}\right)
=\frac{\Lambda}{q}\left(  1-e^{(-\mu-\sqrt{\mu^{2}+2q\sigma^{2}})x/\sigma^{2}%
}\right)  . \label{V en c=0}%
\end{equation}
Let us now state the comparison result for viscosity solutions. The proof is
in the appendix. 
%
%*************
%
%We keep the proof of the Lemma or we say that it is similar o the other paper?
%
%*****************************************

\begin{lemma}
\label{Lema para Unicidad} Assume that (i) $\underline{u}$ is a viscosity
subsolution and $\overline{u}$ is a viscosity supersolution of the HJB
equation (\ref{HJB equation}) for all $x>0$ and for all $c\in\lbrack
0,\overline{c}]$ (ii) $\underline{u}$ and $\overline{u}$ are non-decreasing in
the variable $x$ and in the variable $c$, (iii) $\underline{u}(0,c)=\overline{u}(0,c)=0$ for
$c\in\lbrack0,\overline{c}]$, $\lim_{x\rightarrow\infty}\underline{u}%
(x,c)\leq(\overline{c}+\Lambda)/q\leq\lim_{x\rightarrow\infty}\overline
{u}(x,c)$ and {(iv) $\underline{u}(x,0)\leq\overline{u}(x,0)$ for $x\geq0$}.
Then $\underline{u}\leq\overline{u}$ in $[0,\infty)\times\lbrack0,\overline
{c}).$
\end{lemma}

{\color{black}
Before stating the characterization result, let us describe the concept
of \textit{optimal thresholds} that plays a role below. For a given emission level $c$, the associated optimal threshold represents the surplus level from which it becomes optimal to maintain the emission rate at the (maximal admissible) level  $c$. In particular,
if the optimal threshold associated with $c$ is equal to $0$, then it is
optimal to maintain the emission rate $c$ until depletion of the remaining
carbon budget, i.e., no further reduction of the emission rate occurs before
exhaustion. This interpretation is consistent with the formal definition of
threshold strategies introduced later in Section~6.
} The following characterization theorem is a direct consequence of the previous
lemma, Remark \ref{lim l0/c} and Proposition \ref{Proposicion Viscosidad}.

\begin{theorem}
\label{Caracterizacion Continua}The optimal value function $V$ is the unique
function non-decreasing in $x$ that is a viscosity solution of
(\ref{HJB equation}) in $(0,\infty)\times\lbrack0,\overline{c})$ satisfying
$V(0,c)=0$ and $\lim_{x\rightarrow\infty}$ $V(x,c)=(c+\Lambda)/q$ for
$c\in\lbrack0,\overline{c}).$
\end{theorem}

%From Definition \ref{Optimal Value Function}, Lemma \ref{Lema para Unicidad},
%and Remark \ref{lim l0/c} together with Proposition
%\ref{Proposicion Viscosidad}, we also get the following verification theorem.
%
%\begin{theorem}
%\label{verification result} Consider $S=[0,\overline{c}]$ and consider a
%family of strategies
%\[
%\left\{  C_{x,c}\in\Pi_{x,c}^{S}:(x,c)\in\lbrack0,\infty)\times\lbrack
%0,\overline{c}]\right\}  .
%\]
%If the function $W(x,c):=J(x;C_{x,c})$ {is non-decreasing in the variable $x$}
%and a viscosity supersolution of the HJB equation (\ref{HJB equation}) in
%$(0,\infty)\times\lbrack0,\overline{c})$ with $W(0,c)=0$ and $\lim
%_{x\rightarrow\infty}W(x,c)=$ $(c+\Lambda)/q,$ then $W$ is the optimal value
%function $V$. Also, if for each $k\geq1$ there exists a family of strategies
%$\left\{  C_{x,c}^{k}\in\Pi_{x,c}^{S}:(x,c)\in\lbrack0,\infty)\times
%\lbrack0,\overline{c}]\right\}  $ such that $W(x,c):=\lim_{k\rightarrow\infty
%}J(x;C_{x,c}^{k})$ is a viscosity supersolution of the HJB equation
%(\ref{HJB equation}) in $(0,\infty)\times\lbrack0,\overline{c})$ with
%$W(0,c)=0$ and $\lim_{x\rightarrow\infty}W(x,c)=(c+\Lambda)/q$, then $W$ is
%the optimal value function $V$.
%\end{theorem}

The following proposition establishes conditions under which the current emission level is not lowered anymore, regardless the surplus level.

\begin{proposition}
	\label{Barrera cero}If $\Lambda\leq\sqrt{\mu^{2}+2q\sigma^{2}}$ and
	$\Lambda+\mu>0$, then the optimal threshold is zero for all $c\in
	\lbrack0,\frac{\mu^{2}+2q\sigma^{2}-\Lambda^{2}}{2(\Lambda+\mu)}]$ . If
	$\Lambda+\mu\leq0$ , then the optimal threshold is equal to zero for all
	$c\geq0$.
\end{proposition}

\noindent\textit{Proof.} Consider the value function corresponding to constant
emissions
\[
u(x,c)=W^{c}(x)=\frac{c+\Lambda}{q}\left(  1-e^{\theta_{1}(c)x}\right)  ,
\]
and substitute this function into the HJB equation (\ref{HJB equation}). Since
$\mathcal{L}^{c}(u)(x,c)=0$, in order for $u(x,c)$ to be a solution of (\ref{HJB equation}), we must have $\partial_{c}u(x,c)\geq0$. A direct
computation yields
\[
q\,\partial_{c}u(x,c)=1-e^{\theta_{1}(c)x}\left(  1-\theta_{1}(c)\frac
{c+\Lambda}{\sqrt{(c-\mu)^{2}+2q\sigma^{2}}}\,x\right)  .
\]
To study the sign of this expression, consider
\[
q\,\partial_{cx}u(x,c)=-\theta_{1}(c)e^{\theta_{1}(c)x}\left(1-\frac{c+\Lambda}%
{\sqrt{(c-\mu)^{2}+2q\sigma^{2}}}\left(  1+\theta_{1}(c)x\right)  \right).
\]
Since $\theta_{1}(c)<0$ and $c+\Lambda\geq0$, we have $\partial_{cx}%
u(x,c)\geq0$ if and only if
\[
\frac{c+\Lambda}{\sqrt{(c-\mu)^{2}+2q\sigma^{2}}}\left(  1+\theta
_{1}(c)x\right)  \leq1.
\]
The inequality holds for all $x\geq0$ if and only if%
\[
c+\Lambda\leq\sqrt{(c-\mu)^{2}+2q\sigma^{2}},
\]
which is equivalent to

\[
2c\left(  \Lambda+\mu\right)  +\Lambda^{2}\leq2q\sigma^{2}+\mu^{2}.
\]
Assume first that $\Lambda+\mu>0$. Then the above inequality holds if and only
if
\[
0\leq c\leq\frac{2q\sigma^{2}+\mu^{2}-\Lambda^{2}}{2\left(  \Lambda
	+\mu\right)  }\,.
\]
In particular, we need $2q\sigma^{2}+\mu^{2}-\Lambda^{2}\geq0$, that is
$\Lambda\leq\sqrt{\mu^{2}+2q\sigma^{2}}$ in which case the optimal threshold
is identically zero. If instead $\mu+\Lambda\leq0$ (and hence $\mu^{2}%
\geq\Lambda^{2}$), then%
\[
2c\left(  \Lambda+\mu\right)  \leq0<2q\sigma^{2}+\mu^{2}-\Lambda^{2}.
\]
So the inequality holds for all $c\geq0$, implying that the optimal threshold
is identically zero. Hence, we have the result.\hfill$\blacksquare$

\begin{remark}\label{rem33}\normalfont
 {\color{black}From the proof of Proposition \ref{Barrera cero}, and assuming uniqueness of the
associated threshold levels, one observes that there are no open intervals
of $c$ with zero optimal threshold when
$c>\frac{\mu^{2}+2q\sigma^{2}-\Lambda^{2}}{2(\Lambda+\mu)}$
provided that $\Lambda\leq\sqrt{\mu^{2}+2q\sigma^{2}}$ and
$\Lambda+\mu>0$. Likewise, there are no open intervals with zero optimal
threshold for $c>0$ when
$\Lambda>\sqrt{\mu^{2}+2q\sigma^{2}}$ and $\Lambda+\mu>0$. If, moreover, the optimal thresholds are non-decreasing with respect to
$c$ (as observed in the numerical examples), this implies that, in the case
$\Lambda+\mu>0$, the optimal threshold is strictly positive for
$c>\frac{\mu^{2}+2q\sigma^{2}-\Lambda^{2}}{2(\Lambda+\mu)}$
whenever $\Lambda\leq\sqrt{\mu^{2}+2q\sigma^{2}}$, and strictly positive
for every $c>0$ whenever
$\Lambda>\sqrt{\mu^{2}+2q\sigma^{2}}$.}

\begin{comment}
The optimal threshold being zero for a certain $c$ means that one will
continue to emit at that rate $c$ until the entire carbon budget is depleted. From the proof of Proposition \ref{Barrera cero}, and assuming uniqueness of the
	threshold curve, one observes that there are no open intervals of $c$ with
	zero optimal threshold when $c>\frac{\mu^{2}+2q\sigma^{2}-\Lambda^{2}%
	}{2(\Lambda+\mu)}$ provided that $\Lambda\leq\sqrt{\mu^{2}+2q\sigma^{2}}$ and
	$\Lambda+\mu>0$. Likewise, there are no open intervals with zero optimal
	threshold for $c>0$ when $\Lambda>\sqrt{\mu^{2}+2q\sigma^{2}}$\ and
	$\Lambda+\mu>0$. If, moreover, the optimal threshold curve is non-decreasing
	in $c$ (as we observe in the numerical examples), this implies that, in case $\Lambda+\mu>0$, the
	optimal threshold is strictly positive for $c>\frac{\mu^{2}+2q\sigma
		^{2}-\Lambda^{2}}{2(\Lambda+\mu)}$ whenever $\Lambda\leq\sqrt{\mu
		^{2}+2q\sigma^{2}}$, and the
		optimal threshold is strictly positive for any $c$ $>0$ whenever $\Lambda>\sqrt{\mu
		^{2}+2q\sigma^{2}}$. 
		
\end{comment}

		In other words, if $\Lambda$ is sufficiently large, then under the optimal abatement strategy the emission rate will reach zero at a positive remaining surplus already. That is, the value of not reducing the carbon budget then exceeds the gain from emitting further. This suggests an interpretation of $\Lambda$ as a sort of sustainability parameter that counterbalances the appetite for immediate carbon budget emissions. The limiting value $\sqrt{\mu^{2}+2q\sigma^{2}}$ thus marks the regime in which the sustainability considerations becomes so dominant (relative to emitting) that emissions are halted even when a positive low budget remains. Note that, due to the diffusion properties of the surplus process, the budget may nevertheless be depleted subsequently.

At the same time, these considerations also clarify how the control problem studied in  this paper (almost) degenerates when $\Lambda=0$. Since Proposition \ref{Barrera cero} remains applicable in this case, we conclude that if $\mu\le 0$, the optimal threshold is always zero -- so the initial carbon emission rate is never reduced -- whereas if  $\mu>0$, the emission rate is reduced only as long as $c>(\mu^2+2q\sigma^2)/(2\mu)$. See also the numerical example in Section 7 for an illustration.% \par\noindent
\hfill$\diamond$
\end{remark}

\begin{remark}
\label{Remark sigma zero}\normalfont Let us consider the limit case where $\sigma=0$ and
$\mu\geq0$. If $c\in\lbrack0,\mu]$, then the surplus is never depleted. So,
it is straightforward to verify that the optimal threshold in this case is
zero. Therefore, the corresponding optimal value function, which results from
emitting at the constant rate $c$ indefinitely, is given by
\[
V(x,c)=\int_{0}^{\infty}(\Lambda+c)e^{-qs}=\left(  c+\Lambda\right)  /q.
\]
Now, consider the case $c>\mu$ and $x>0$. An admissible strategy in this
setting is to maintain constant emissions at the maximum admissible level $c$
while the surplus is positive, i.e. for $0\leq t<$ $T:=x/(c-\mu)$. Once the
surplus hits zero at time $T$, the emission rate is reduced to the (maximum possible) level $\mu\geq0$ which can be sustained indefinitely. The value
function for this strategy is
\begin{align*}
W(x,c)  &  =\int_{0}^{x/(c-\mu)}\left(  c+\Lambda\right)  e^{-qs}%
ds+e^{-qx/(c-\mu)}\frac{\mu+\Lambda}{q}\\
%&  =\frac{c+\Lambda}{q}\left(  1-e^{-qx/(c-\mu)}\right)  +\frac{\mu+\Lambda
%}{q}e^{-qx/(c-\mu)}\qquad (\ast)\\
&  =\frac{c+\Lambda}{q}-\frac{c-\mu}{q}e^{-qx/(c-\mu)}\,.
\end{align*}
Due to the discount factor $q>0$ and the fact that the emission level can
be reduced to $c_{0}=\mu\geq0$ at zero surplus (the surplus can then not
become negative unlike in the Brownian setting), it is clear that this is the
optimal strategy and so $W$ is the optimal value function. We now relate this
result to the HJB framework. While a full proof of the HJB approach is
omitted in this simplified setting, we can infer that the corresponding HJB
equation is
\[
\max\{\overline{\mathcal{L}}^{c}(W)(x,c),-\partial_{c}W(x,c)\}=0
\]
for $(x,c)\in\lbrack0,\infty)\times(\mu,\overline{c}],$ where
\[
\overline{\mathcal{L}}^{c}(W)(x,c):=(\mu-c)\partial_{x}%
W(x,c)-qW(x,c)+c+\Lambda=0.
\]
This corresponds to put $\sigma=0$ in $\mathcal{L}^{c}$ and the boundary
condition $W(0,c)=\left(  \mu+\Lambda\right)  /q>0$. The latter reflects the fact
that, even at zero surplus, it is possible to emit $c=\mu\geq0$ indefinitely.
$W$ satisfies the associated first-order HJB equation. To see this, it is
immediate to show that $W$ is a solution of $\overline{\mathcal{L}}%
^{c}(W)(x,c)=0$. Additionally, differentiating $W$ with respect to $c$
yields
\[
\partial_{c}W(x,c)=\frac{1}{q}\left(  1-e^{-qx/(c-\mu)}\left(  1+\frac
{qx}{c-\mu}\right)  \right)
\]
for $\left(  x,c\right)  \in(\mu,\overline{c}]$. Setting $a=qx/(c-\mu
)\geq0$, we get the inequality $1\geq e^{-a}\left(  1+a\right)$, which implies
$\partial_{c}W(x,c)\geq0.$ Hence, the usual verification condition holds.\hfill $\diamond$
\end{remark}

\section{Hamilton-Jacobi-Bellman equations for finite sets
\label{Hamilton-Jacobi-Bellman equations for finite sets}}

Let us now \ restrict to the following finite set of possible emission rates:%
\[
S=\left\{  c_{0},c_{1},c_{2},....,c_{n}\right\},
\]
where $0=c_{0}<c_{1}<c_{2}<....<c_{n}=\overline{c}$. Note  that $V^{S}%
(x,c_{i})=V^{\left\{  0,c_{1},....,c_{i}\right\}  }(x,c_{i})$, i.e., it
depends only on the emission rates up to $c_{i}$ and does not involve
$c_{i+1},...,c_{n}$. To simplify the notation, we define the optimal value
function within the finite set $S$%

\begin{equation}
V^{c_{i}}(x):=V^{S}(x,c_{i}), \label{optimal value function Discreta}%
\end{equation}
which represents the optimal value function corresponding to initial emission
level $c_{i}$. We then have the following inequalities:
\[
V^{c_{i}}(x)\geq V^{c_{i-1}}(x)\geq...\geq V^{c_{0}}(x)\geq0,
\]
where
\[
V^{0}(x)=V^{\{0\}}(x,0)=\frac{\Lambda}{q}\left(  1-e^{\theta_{1}(0)x}\right)
.
\]

Assuming $C^{2}([0,\infty))$-regularity for $V^{c_{i}}$, we can heuristically
derive the HJB equation associated to the discrete optimal
value function (\ref{optimal value function Discreta}):
\begin{equation}
\max\left\{  \mathcal{L}^{c_{i}}(v(x)),V^{c_{i-1}}(x)-V^{c_{i}}(x)\right\}
=0\text{ for }x\geq0\text{ and }i=1,...,n\text{,}
\label{HJB equation discreta}%
\end{equation}
with $V^{c_{i}}(0)=0$ and $\lim_{x\rightarrow\infty}$ $V^{c_{i}}%
(x)=(\Lambda+c_{i})/q$. Let us define
\[
v^{0}=V^{0}=\frac{\Lambda}{q}\left(  1-e^{\theta_{1}(0)x}\right)
\]
and the system of ODE's%
\begin{equation}
\max\left\{  \mathcal{L}^{c_{i}}(v^{i}(x)),v^{i-1}(x)-v^{i}(x)\right\}
=0\text{ for }x\geq0\text{ and }i=1,...,n\text{,} \label{Discrete System ODE}%
\end{equation}
with $v^{i}(0)=0$ and $\lim_{x\rightarrow\infty}$ $v^{i}(x)=(\Lambda+c_{i})/q$.\\

Let us now show that $V^{c_{i}}$ is the unique solution in the viscosity sense
to the ODE system (\ref{Discrete System ODE}). For this purpose, let
us introduce first the definition of a viscosity solution in the
one-dimensional case.

\begin{definition}
\label{Viscosidad Discreta}(a) A locally Lipschitz function $\overline
{u}:[0,\infty)\rightarrow{\mathbb{R}}$\ is a viscosity supersolution of
(\ref{HJB equation discreta})\ at $x\in(0,\infty)$ if any twice continuously
differentiable function $\varphi:[0,\infty)\rightarrow{\mathbb{R}}\ $with
$\varphi(x)=\overline{u}(x)$, such that $\overline{u}-\varphi$\ reaches the
minimum at $x$, satisfies
\[
\max\left\{  \mathcal{L}^{c_{i}}(\varphi)(x),V^{c_{i-1}}(x)-\varphi
(x)\right\}  \leq0.\
\]
The function $\varphi$ is called a \textbf{test function for supersolution} at
$x$.

(b) A locally Lipschitz function $\underline{u}:$ $[0,\infty)\rightarrow
{\mathbb{R}}\ $\ is a viscosity subsolution\ of (\ref{HJB equation discreta}%
) at $x\in(0,\infty)$ if any twice continuously differentiable function
$\psi:[0,\infty)\rightarrow{\mathbb{R}}$ with $\psi(x)=\underline{u}(x)$, such
that $\underline{u}-\psi$ reaches the maximum at $x$,  satisfies
\[
\max\left\{  \mathcal{L}^{c_{i}}(\psi)(x),V^{c_{i-1}}(x)-\psi(x)\right\}
\geq0\text{.}%
\]
The function $\psi$ is called a \textbf{test function for subsolution} at $x$.

(c) A function $u:[0,\infty)\rightarrow{\mathbb{R}}$ which is both a
supersolution and subsolution at $x\in\lbrack0,\infty)$ is called a viscosity
solution of (\ref{HJB equation discreta})\ at $x$.
\end{definition}

The following characterization theorem is the discrete analogue of Theorem
\ref{Caracterizacion Continua}. The proof is omitted, as it follows similar
arguments to those in the continuous case but is technically simpler.

\begin{theorem}
\label{Caracterizacion Discreta}The optimal value function $V^{c_{i}}(x)$ for
$1\leq i<n$ is the unique viscosity solution of the associated HJB equation
(\ref{HJB equation discreta}) with boundary condition $V^{c_{i}}(0)=0$ and
$\lim_{x\rightarrow\infty}V^{c_{i}}(x)=(\Lambda+c_{i})/q.$
\end{theorem}

We also have the following alternative characterization theorem.

\begin{theorem}
\label{Menor supersolucion discreta}The optimal value function $V^{c_{i}}(x)$
for $0\leq i<n$ is the smallest viscosity supersolution of the associated HJB
equation (\ref{HJB equation discreta}) with boundary condition $0$ at $x=0$
and limit greater than or equal to $(\Lambda+c_{i})/q\ $as $x$ goes to infinity.
\end{theorem}

Since for $i\geq1,$ the optimal value function $V^{c_{i}}$ is a viscosity
solution of (\ref{HJB equation discreta}), there are values of $x$ where
$V^{c_{i}}(x)=V^{c_{i-1}}(x)$ and values of $x$ where $\mathcal{L}^{c_{i}%
}(V^{c_{i}})(x)=0$. So for any $i\geq1$, we can { partition} $(0,\infty)$
{into the closed} set $\mathcal{D}_{i}=\{x:V^{c_{i}}(x)=V^{c_{i-1}}(x)\}$ and
the open set $\mathcal{E}_{i}=\{x:V^{c_{i}}(x)>V^{c_{i-1}}(x)\}.$ Moreover,
$\mathcal{L}^{c_{i}}(V^{c_{i}})(x)=0$ in $\mathcal{E}_{i}$ and the optimal
strategy is to emit at rate $c_{i}$ when the current
surplus is in $\mathcal{E}_{i}$ and to decrease the emission rate when the
current surplus is in $\mathcal{D}_{i}$.

\section{Convergence of the optimal value functions from the discrete to the
continuous case \label{Seccion Convergencia}}

In this section, we prove that the optimal value functions corresponding to
the (ratcheting-down) finite set of possible carbon emission rates, as
defined in the previous section, converge to the optimal value function of the
continuous case as the mesh size of the finite sets approaches zero. This is
achieved by considering a sequence of nested meshes. \\
%
%In the next section, we
%will demonstrate that the optimal strategies corresponding to the finite set
%of $CO_{2}$.

Consider, for $n\geq0$, a sequence of sets $\mathcal{S}^{n}$ (each with
$k_{n}$ elements) of the form
\[
\mathcal{S}^{n}=\left\{  c_{0}^{n}=0<c_{k_{1}}^{n}<\cdots<c_{k_{n}}%
^{n}=\overline{c}\right\}
\]
satisfying the conditions $\mathcal{S}^{0}=\left\{  0,\overline{c}\right\}  $,
$\mathcal{S}^{n}\subset\mathcal{S}^{n+1}$ and mesh-size $\delta(\mathcal{S}%
^{n}):=\max_{k=1,k_{n}}\left(  c_{k}^{n}-c_{k-1}^{n}\right)  \searrow0$ as $n$
goes to infinity. We extend the definition of $V^{\mathcal{S}^{n}}$ to a
function {$V^{n}:[0,\infty)\times\lbrack 0,\overline{c}%
]\rightarrow{\mathbb{R}}$} as follows:%
\begin{equation}
V^{n}(x,c)=V^{\mathcal{S}^{n}}(x,\widetilde{c}^{n}), \label{Definicion Vn}%
\end{equation}
where
\begin{equation}
\widetilde{c}^{n}=\max\{c_{i}^{n}\in\mathcal{S}^{n}:c_{i}^{n}\leq c\}.
\label{Deminicioncruliton}%
\end{equation}
We will prove that $\lim_{n\rightarrow\infty}V^{n}(x,c)=V^{[0,\overline{c}%
]}(x,c)$ for any $(x,c)\in\lbrack0,\infty)\times\lbrack0,\overline{c}]$ and we
will study the uniform convergence of this limit. Since $\mathcal{S}%
^{n}\subset\mathcal{S}^{n+1}$, it follows that $\widetilde{c}^{n+1}\leq$
$\widetilde{c}^{n}\in\mathcal{S}^{n}$ for each $c\in\lbrack 0,\overline{c}]$. Then, by monotonicity of $V^{\mathcal{S}^{n+1}}$ with
respect to its second variable,

\[
V^{[0,\overline{c}]}(x,c)\geq V^{n+1}(x,c)=V^{\mathcal{S}^{n+1}}%
(x,\widetilde{c}^{n+1})\geq V^{\mathcal{S}^{n+1}}(x,\widetilde{c}^{n})\geq
V^{\mathcal{S}^{n}}(x,\widetilde{c}^{n})=V^{n}(x,c).
\]
Therefore, the pointwise limit exists and we can define the limit function as
\begin{equation}
\overline{V}(x,c):=\lim_{n\rightarrow\infty}V^{n}(x,c).
\label{Definicion Vbarra}%
\end{equation}
Later on, we will show that $\overline{V}=V^{[0,\overline{c}]}$. Note that
$\overline{V}(x,c)$ is non-increasing in $c$, satisfies $\overline
{V}(x,\overline{c})=V(x,\overline{c})$, and is non-decreasing in $x$, with the
asymptotic behavior $\lim_{x\rightarrow\infty}$ $\overline{V}(x,c)=(\overline
{c}+\Lambda)/q$.\\

Using the same arguments as in Proposition 6.1 of Albrecher et al.\ \cite{AAM}, we obtain the following result.

\begin{proposition}
\label{Proposicion Convergencia Uniforme VN a V}The sequence $V^{n}$ converges
uniformly to $\overline{V}.$
\end{proposition}

{ }

Note that for any {$(x,c)\in\lbrack0,\infty)\times\lbrack0,\overline{c}]$}, we
have that $V^{n}(x,c)=V^{\mathcal{S}^{n}}(x,\widetilde{c}^{n})$ is a value
function corresponding to an admissible strategy in $\Pi_{x,\widetilde
{c}^{n}}^{\mathcal{S}^{n}}\subset\Pi_{x,c}^{[0,\overline{c}]}$. Hence
$\overline{V}(x,c)=\lim_{n\rightarrow\infty}V^{n}(x,c)$ is itself a limit of
value functions of admissible strategies in $\Pi_{x,c}^{[0,\overline{c}]}$.
Moreover, by Proposition \ref{Proposition Global Lipschitz zone},%
\[
0\leq V^{n}(x_{2},c_{1})-V^{n}(x_{1},c_{2})\leq K\left[  \left(  x_{2}%
-x_{1}\right)  +\left(  c_{2}-c_{1}\right)  \right]
\]
for all $n$, with a constant $K$ independent on $n$. Since $V^{n}$ converges
uniformly to $\overline{V}$, it follows that $\overline{V}$ is Lipschitz with
the same constant $K.$\\

With this result, we are now in a position to state the main result of this
section. We omit the proof, as it closely follows the one given in Theorem 4.2
of Albrecher et al.\ \cite{AAM1}.

\begin{theorem}
\label{Vbarra es V} The function $\overline{V}$ defined in
(\ref{Definicion Vbarra}) is the optimal value function $V^{[0,\overline{c}]}$.
\end{theorem}

\section{Optimal strategies for finite sets \label{Seccion optima Discreta}}

Let us once again consider a finite set of possible 
emissions rates:%
\[
S=\left\{  c_{0},c_{1},c_{2},....,c_{n}\right\}  ,
\]
where $0=c_{0}<c_{1}<c_{2}<....<c_{n}=\overline{c}$. { We first look for the
following particular strategies, which we call multi-threshold strategies.
These are defined as follows: }

\begin{itemize}
\item { }${v}^{0}${$(x):=V^{0}(x)=$ }$\frac{\Lambda}{q}\left(  1-e^{\theta
_{1}(0)x}\right)  $

\item {For each }$i\geq1$ and {thresholds }$z(c_{i})\geq0$, {the value
function $v^{c_{i}}(x)$ satisfies $\mathcal{L}^{c_{i}}(v^{c_{i}})(x)=0$ for
$x\in(z(c_{i}),\infty)$ with }$\lim_{x\rightarrow\infty}W^{z}(x,c_{i}%
)=\frac{c_{i}+\Lambda}{q}$ and $v^{c_{i}}(x)=v^{c_{i-1}}(x)$ for $x\in
\lbrack0,z(c_{i})]$.
\end{itemize}

We will show in this section that the optimal discrete value functions
$V^{c_{i}}$ are indeed of this form. As a result, the optimal value function in
the continuous control setting can be seen as the limit of value functions
associated with multi-threshold strategies.\\
{\color{black} Since no further reduction of the emission rate is
possible once the rate $0$ is reached, the threshold function is defined
only on the strictly positive emission rates. For this reason, we introduce
$\widetilde{S}=\{c_{1},c_{2},\dots,c_{n}\}$ separately and consider a
function $z:\widetilde{S}\rightarrow[0,\infty)$.}

\begin{comment}
More precisely, let $\widetilde{S}=\left\{  c_{1},c_{2},....,c_{n}\right\}  $
and consider a function $z:\widetilde{S}\rightarrow\lbrack0,\infty)$.
\end{comment}

We then
define a \textit{threshold strategy} (which depends on both the current surplus $x$ and the emission rate $c_{i}\in S$), recursively as a stationary strategy%

\begin{equation}
\mathbf{\pi}^{z}=(C_{x,c_{i}})_{(x,c_{i})\in\lbrack0,\infty)\times {\color{black}\widetilde{S}}}\text{
where }C_{x,c_{i}}\in\Pi_{x,c_{i}}^{{\color{black}\widetilde{S}}} \label{Pizzeta}%
\end{equation}
as follows:

\begin{itemize}
\item If $i=0$ (i.e. no carbon emission), then $(C_{x,{\color{black}c_{0}}})_{t}=0$.

\item If $1\leq i\leq n$ and $x\leq z(c_{i})$ with $z(c_{i})\geq z(c_{i-1})$,
follow $C_{x,c_{i-1}}\in\Pi_{x,c_{i-1}}^{{\color{black}\widetilde{S}}}$.

\item If $1\leq i\leq n$ and $x>z(c_{i})$ emit with rate $c_{i}$ as
long as the surplus exceeds $z(c_{i})$; once the current surplus reaches
$z(c_{i})$, switch to $C_{x,c_{i-1}}\in\Pi_{x,c_{i-1}}^{{\color{black}\widetilde{S}}}$. More precisely,
{
\[
(C_{x,c_{i}})_{t}=c_{i}I_{t<\widehat{\tau}}+(C_{z(c_{i}),c_{i-1}}%
)_{t~}I_{\widehat{\tau}\leq t<\tau},
\]
} where $\widehat{\tau}$ is the first hitting time of the surplus process to
the level $z(c_{i})$ and $\tau$ is the depletion time$.$
\end{itemize}

We refer to $z(c_{i})$ as the \textit{threshold at emission rate level}
$c_{i}$ and the function $z:\widetilde{S}\rightarrow\lbrack0,\infty)$ as the
\textit{threshold function. }The expected payoff of the multi-threshold
strategy $\pi^{z}$ is given by
\begin{equation}
W^{z}(x,c_{i}):=J(x;C_{x,c_{i}}). \label{Wz discreta}%
\end{equation}
Note that $W^{z}(x,c_{i})$ only depends on the threshold values $z(c_{k})$ for
$1\leq k\leq i,$ that $W^{z}(0,c_{i})=0$ and that $W^{z}(x,c_{i})=V^{0}(x)$
for $0\leq x\leq\min\{z(c_{k}):1\leq k\leq i\}.$

We next obtain a recursive formula for $W^{z}$.

\begin{proposition}
\label{Formula recursiva de Wz discreta} We have the following recursive
formula for $W^{z}$: $W^{z}(x,0)=\Lambda\left(  1-e^{\theta_{1}(0)x}\right)
/q$, and for $1\leq i\leq n,$
\[
W^{z}(x,c_{i})=\left\{
\begin{array}
[c]{lll}%
W^{z}(x,c_{i-1}) & \text{if} & x\leq z(c_{i})\\
\frac{c_{i}+\Lambda}{q}\left(  1-a^{z}(c_{i})e^{\theta_{1}(c_{i})x}\right)
\  & \text{if} & x>z(c_{i}),
\end{array}
\right.
\]
where%
\[
a^{z}(c_{i}):=\left(  1-\frac{q}{c_{i}+\Lambda}W^{z}(z(c_{i}),c_{i-1})\right)
e^{-\theta_{1}(c_{i})z(c_{i})}\text{ and }e^{\theta_{1}(c_{i})z(c_{i})}%
>a^{z}(c_{i})>0.
\]

\end{proposition}

\noindent\textit{Proof.} By construction, the strategy $\mathbf{\pi}^{z}$ emits at rate
$c_{i}$ when the surplus exceeds $z(c_{i}).$ Hence, 
$\mathcal{L}^{c_{i}}(W^{z})(x,c_{i})=0$ for $x\in(z(c_{i}),\infty)$. Since
$\lim_{x\rightarrow\infty}W^{z}(x,c_{i})=(c_{i}+\Lambda)/q$ and the
emission strategy switches to emit $c_{i-1}$ at the threshold $z(c_{i})$, we
have $W^{z}(z(c_{i}),c_{i})=W^{z}(z(c_{i}),c_{i-1})$. Also, $W^{z}%
(x,c_{i-1})<(c_{i}+\Lambda)/q$, so we get the result. \hfill
$\blacksquare$\newline

\bigskip

Now, we aim to maximize $W^{z}(x,c_{i})$ over all possible multi-threshold
functions\textit{\ }$z:\widetilde{S}\rightarrow\lbrack0,\infty)$. We denote by
$z^{\ast}:\widetilde{S}\rightarrow\lbrack0,\infty)$ the optimal
multi-threshold function, which can equivalently be interpreted as the one
that minimizes $a^{z^{\ast}(c_{i})}(c_{i})$ for each $1\leq
i\leq n$. From Proposition \ref{Barrera cero}, if $\Lambda<\sqrt{\mu
^{2}+2q\sigma^{2}}$, then $z^{\ast}(c_{i})=0$ for all $c_{i}\in\lbrack
0,\frac{\mu^{2}+2q\sigma^{2}-\Lambda^{2}}{2(\Lambda+\mu)}]$ and $z^{\ast
}(c_{i})>0$ otherwise. Therefore, from now on we consider only the case
$c_{i}>\frac{\mu^{2}+2q\sigma^{2}-\Lambda^{2}}{2(\Lambda+\mu)}$ if this value
is positive. Note that, as a first step, we are maximizing the discounted
expected emissions only among multi-threshold strategies, not among all
admissible strategies, which could, in principle, have a more complex
structure. {Later in this section (Theorem \ref{Optima Discreta Threshold}),
we will show that the resulting value function coincides with the optimal
discrete value function $V^{S}(x,c_{i})$.}

Since the initial function $W^{z}(x,0)$ in the recursive procedure is known,
we can interpret the optimization problem in two different ways.

\begin{enumerate}
\item First Approach. Recursive One-Dimensional Optimization:

We solve a sequence of $n$ one-dimensional optimization problems obtaining the
minimum of $a^{z}(c_{i})$. Suppose that $W^{z^{\ast}}(x,c_{k})$ and $z^{\ast
}(c_{k})$ are known for $k=1,\ldots,i-1$. Then, from the recursive formula
(Proposition \ref{Formula recursiva de Wz discreta}), we can
compute $W^{z^{\ast}}(x,c_{i})$ and $z^{\ast}(c_{i})$ as follows. Define the
continuous function $G_{i}:[0,\infty)\rightarrow\mathbb{R}$ as%
\begin{equation}
G_{i}(y):=\left(  1-\frac{q}{c_{i}+\Lambda}W^{z^{\ast}}(y,c_{i-1})\right)
e^{-\theta_{1}(c_{i})y}\text{ .} \label{Definicion Gi (1)}%
\end{equation}
We have $G_{i}(0)=1$ and since
\[
0<\lim_{y\rightarrow\infty}W^{z^{\ast}}(y,c_{i-1})<\frac{c_{i-1}+\Lambda}{q}%
\]
and $\theta_{1}(c_{i})<0$, we have $0<G_{i}(y)$ and $\lim_{y\rightarrow\infty
}G_{i}(y)=\infty$. As $G_{i}$ is continuous, it attains its minimum
in $[0,\infty)$. We define
\begin{equation}
z^{\ast}(c_{i}):=\min\left(  \arg\min_{y\in\lbrack0,\infty)}G_{i}(y)\right)
,\text{ }a^{\ast}(c_{i})=G_{i}(z^{\ast}(c_{i})). \label{First Approach}%
\end{equation}
The function $W^{z^{\ast}}(\cdot,c_{i})$ satisfies $\mathcal{L}^{c_{i}%
}(W^{z^{\ast}})(x,c_{i})=0$ for $x>z^{\ast}(c_{i})$ and $W^{z^{\ast}}%
(x,c_{i})=$ $W^{z^{\ast}}(x,c_{i-1})$ for $x\in\lbrack0,z^{\ast}(c_{i})].$

\item Second Approach: Sequence of Obstacle Problems.

This approach interprets the problem as a recursive sequence of
one-dimensional obstacle problems. Assume that $W^{z^{\ast}}(x,c_{k})$ and the
optimal thresholds $z^{\ast}(c_{k})$ are known for $k=1,\ldots,i-1$. To find
$W^{z^{\ast}}(x,c_{i})$ and $z^{\ast}(c_{i})$, consider the smallest solution
$U^{\ast}$ of the differential equation $\mathcal{L}^{c_{i}}(U)=0$ in
$[0,\infty)$ with boundary condition $\lim_{x\rightarrow\infty}U(x)=\frac
{c_{i}+\Lambda}{q}$ such that $U^{\ast}(\cdot)\geq W^{z^{\ast}}(\cdot
,c_{i-1})$. We define:%
\begin{equation}
z^{\ast}(c_{i})=\left\{
\begin{array}
[c]{ll}%
0 & \text{if }U^{\ast}(\cdot)>W^{z^{\ast}}(\cdot,c_{i-1})\ \text{in }%
(0,\infty)\\
\sup\{y>0:U^{\ast}(y)=W^{z^{\ast}}(y,c_{i-1})\} & \text{otherwise.}%
\end{array}
\right.  \label{Second Approach}%
\end{equation}

In other words, $z^{\ast}(c_{i})$ is the last point at which $U^{\ast}$ and
$W^{z^{\ast}}(\cdot,c_{i+1})$ coincide. If they only coincide at $y=0$, then
$z^{\ast}(c_{i})=0$. We then have that $W^{z^{\ast}}(x,c_{i})=U^{\ast}(x)$ for
$x>z^{\ast}(c_{i})$ and $W^{z^{\ast}}(x,c_{i})=W^{z^{\ast}}(x,c_{i-1})$ for
$x\leq z^{\ast}(c_{i})$. To show that $U^{\ast}$ exists, note that by
(\ref{Bounded solutions}), the solutions $U$ of the
differential equation $\mathcal{L}^{c_{i}}(U)=0$ in $[0,\infty)$ with boundary
condition $\lim_{x\rightarrow\infty}U(x)=(c_{i}+\Lambda)/q$ are of the
form
\[
U_{a}(x)=\frac{c_{i}+\Lambda}{q}\left(  1-ae^{\theta_{1}(c_{i})x}\right)  .
\]
So, $U^{\ast}=U_{a^{\ast}(c_{i})}$ where $a^{\ast}(c_{i})$ is defined in
(\ref{First Approach}).
\end{enumerate}

\begin{remark}
\normalfont\label{Remark Obstaculo} If $z^{\ast}(c_{i})>0,$ $U_{a_{i}^{\ast}%
}(z^{\ast}(c_{i}))=W^{z^{\ast}}(z^{\ast}(c_{i}),c_{i-1})$, $U_{a_{i}^{\ast}%
}(x)\geq W^{z^{\ast}}(x,c_{i-1})$ for $x\geq0$ and $U_{a_{i}^{\ast}%
}(x)>W^{z^{\ast}}(x,c_{i+1})$ for $x\in(z^{\ast}(c_{i}),\infty)$. Note that we
can show by a recursive argument that $W^{z^{\ast}}(x,c_{i})$ is infinitely
continuously differentiable at all $x\in\lbrack0,\infty)\setminus\{z^{\ast
}(c_{k}):k=i,\ldots,n\}$ and continuously differentiable at the points $z^{\ast
}(c_{k})\ $for $k=i,\ldots,n.$ Indeed, $U_{a_{i}^{\ast}}$ and $W^{z^{\ast}}%
(\cdot,0)$ are infinitely continuously differentiable and $U_{a_{i}^{\ast}%
}^{\prime}(z^{\ast}(c_{i}))-\partial_{x}W^{z^{\ast}}(z^{\ast}(c_{i}%
),c_{i-1})=0$ because $U_{a_{i}^{\ast}}(\cdot)-W^{z^{\ast}}(\cdot,c_{i-1})$
reaches the minimum at $z^{\ast}(c_{i})$. Moreover, since $W^{z^{\ast}%
}(x,c_{i})=W^{z^{\ast}}(x,c_{i-1})I_{\{x<z^{\ast}(c_{i})\}}+U_{a_{i}^{\ast}%
}(x)I_{\{x\geq z^{\ast}(c_{i})\}}$,
\[
\partial_{xx}W^{z^{\ast}}(z^{\ast}(c_{i})^{+},c_{i})-\partial_{xx}W^{z^{\ast}%
}(z^{\ast}(c_{i})^{-},c_{i})=U_{a_{i}^{\ast}}^{\prime\prime}(z^{\ast}%
(c_{i}))-\partial_{xx}W^{z^{\ast}}(z^{\ast}(c_{i})^{-},c_{i-1})\geq0.
\]
\hfill $\diamond$
\end{remark}

\begin{remark}\normalfont 
\label{Remark Concavity}The function $U_{0}(x)=(c_{i}+\Lambda)/q$ is a constant
function. For $a>0$, the function
\[
U_{a}(x)=\frac{c_{i}+\Lambda}{q}\left(  1-ae^{\theta_{1}(c_{i})x}\right)  .
\]
is strictly increasing and concave, with 
\[\partial_{x}U_{a}(x)=-\theta_{1}(c_{i})\frac{c_{i}+\Lambda}{q}ae^{\theta
_{1}(c_{i})x}>0,\quad \partial_{xx}U_{a}(x)=-\theta_{1}^{2}(c_{i}%
)\frac{c_{i}+\Lambda}{q}ae^{\theta_{1}(c_{i})x}<0,
\]
and it is bounded above by $U_{0}(x).$\hfill $\diamond$
\end{remark}

In the next theorem, we show that there exists an optimal strategy and it is
of threshold type. The proof is in the appendix.

\begin{theorem}
\label{Optima Discreta Threshold} If $z^{\ast}\ $is the optimal threshold
function, then $W^{z^{\ast}}(x,c_{i})$ is the optimal function $V^{c_{i}}(x)$
defined in (\ref{Optimal Value Function}) for $i=1,...,n$.
\end{theorem}

\begin{remark}
\label{Formula threshold discreta} By Remark \ref{Remark Obstaculo}, the
function $G_{i}$ defined in (\ref{Definicion Gi (1)})$\ $ is continuously
differentiable. If its minimum $z^{\ast}(c_{i})$ is positive, the first-order
condition $G_{i}^{\prime}(z^{\ast}(c_{i}))=0$ implies that $z^{\ast}(c_{i}%
)$ satisfies the implicit equation

%\begin{proposition}%
\[
\theta_{1}(c_{i})W^{z^{\ast}}(y,c_{i-1})-\partial_{x}W^{z^{\ast}}%
(y,c_{i-1})=\theta_{1}(c_{i})\frac{c_{i}+\Lambda}{q}%
\]
for $i=1,\ldots,n-1.$ \hfill $\diamond$
\end{remark}
%\end{proposition}

\begin{remark}
\normalfont\label{Pizzeta extendido} Given $z:\widetilde{S}\rightarrow
\lbrack0,\infty)$, we have defined in (\ref{Pizzeta}) a threshold strategy
$\mathbf{\pi}^{z}=(C_{x,c_{i}})_{(x,c_{i})\in\lbrack0,\infty)\times S}$, where
$C_{x,c_{i}}\in\Pi_{x,c_{i}}^{S}$ for $i=1,\ldots,n$. We can extend this
threshold strategy to
\begin{equation}
\widetilde{\mathbf{\pi}}^{z}=(C_{x,c})_{(x,c)\in\lbrack0,\infty)\times
\lbrack0,c_{n}]}~\text{where~}C_{x,c}\in\Pi_{x,c}^{S} \label{Pizzeta Rulo}%
\end{equation}
as follows:

\begin{itemize}
\item If $c\in(c_{i},c_{i+1})$ and $x>z(c_{i})$, emit with rate $c$ while the
current surplus is above $z(c_{i})$. If the current surplus reaches
$z(c_{i})$, follow $C_{z(c_{i}),c_{i}}\in\Pi_{x,c_{i}}^{S}$.

\item If $c\in(c_{i},c_{i+1})$ and $x\leq z(c_{i})$ , follow $C_{x,c_{i}}%
\in\Pi_{x,c_{i}}^{S}.$ More precisely, if $(x,c)\in\lbrack z(c_{i}),\infty)\times(c_{i},c_{i+1}),$
then $C_{x,c}\in\Pi_{x,c}^{S}$ is defined as $\left(  C_{x,c}\right)  _{t}=c$
and so $X_{t}^{C_{x,c}}=X_{t}-ct$ for $t<\tau_{i}$ where
\[
\tau_{i}:=\min\{s:X_{t}^{C_{x,c}}=z(c_{i})\},
\]
and $\left(  C_{x,c}\right)  _{t}=\left(  C_{z(c_{i}),c_{i}}\right)
_{t-\tau_{i}}\in\Pi_{z(c_{i}),c_{i}}^{S}$ for $t\geq\tau_{i}$. Finally,
$C_{x,c}=C_{x,c_{i}}$ $\in\Pi_{x,c_{i}}^{S}$ for $(x,c)\in\lbrack
0,z(c_{i})]\times(c_{i},c_{i+1}).$
\end{itemize}

The value function of the stationary strategy $\widetilde{\mathbf{\pi}}^{z}$
is defined as
\begin{equation}
J^{\widetilde{\mathbf{\pi}}^{z}}(x,c):=J(x;C_{x,c}):[0,\infty)\times\lbrack
c_{0},c_{n}]\rightarrow\lbrack0,\infty){.} \label{W Pizzeta Rulo}%
\end{equation}
\hfill $\diamond$
\end{remark}

\section{Numerical Illustrations \label{Numerical examples}}

In this section we present examples in which we approximate the optimal value
function by a multi-threshold strategy considering a discrete set of possible
emission rates. For a given $n$, define the mesh-size as $\Delta
c={\overline{c}}/{n}$ and consider the finite set%

\[
S^{n}=\left\{  0,\Delta c,2\Delta c,3\Delta c,....,\overline{c}\right\}.
\]

(1) We begin by defining
\[
V^{0}(x)=\frac{\Lambda}{q}\left(  1-e^{\theta_{1}(0)x}\right),
\]
which is the solution to the equation $\mathcal{L}^{0}(W)=0$ with limit
boundary conditions $\lim_{x\rightarrow\infty}W(x)=\Lambda/q$. Note that
$V^{0}$ is not zero due to the positive reward $\Lambda$.\\

(2) Recursive construction:

To compute $V^{k\Delta c},$ we consider value functions of strategies that
emit at a constant rate $c_{k}=k\Delta c$ when $x\geq x_{k}$ and
switch to the lower rate value function $V^{(k-1)\Delta c}$ when $0\leq
x<x_{k}$. To obtain this value function, we consider the solutions of equation
$\mathcal{L}^{k\Delta c}(W_{1})=0$ on $(x_{k},\infty]$ with boundary condition
at infinity $\lim_{x\rightarrow\infty}W_{k}(x)=(k\Delta c+\Lambda)/q.$ The
general solution is given by
\[
W_{k}(x)=\frac{k\Delta c+\Lambda}{q}+a_{k}e^{\theta_{1}(k\Delta c)x}.
\]
We then determine the constant $a_{k}$ by matching this function continuously
to $V^{(k-1)\Delta c}(x)$ at the threshold point $x_{k}$. Finally, we optimize
over all possible switching points $x_{k}$ to obtain the optimal threshold
$z^{\ast}(c_{k})=x_{k}^{\ast.}.$ It follows that $V^{k\Delta c}(x)$ is the
optimal value function corresponding to the optimal multi-threshold strategy
described above.

In each of the examples, we display $V^{S}(x,\overline{c})$ as a function of initial carbon budget $x$. From the results of the
previous sections, we know that this function converges to the optimal value function of the 
continuous case as $n\rightarrow\infty$ (we choose $n=500$ in each of the illustrations). We also depict the set $\left\{
(z^{\ast}(0),0),(z^{\ast}(c_{1}),c_{1}),...,(z^{\ast
}(\overline{c}),\overline{c})\right\}  $ that corresponds to the optimal threshold points.  {\color{black}These points are then used to approximate the optimal strategy in the continuous-control setting, which is characterized by a free-boundary curve. More precisely, to the right of this curve, it is optimal to continue emitting carbon at the highest admissible rate, namely the current emission level $c$. In contrast, to the left of the curve, the optimal policy is to reduce the emission rate immediately, thereby moving vertically downward in the \((x,c)\)-plane until the state reaches the curve.
}\\
\begin{comment}  
These points
are taken as an approximation of the optimal strategy---hence of a curve---in
the continuous case. Roughly speaking, to the right of this curve, the optimal
strategy is to emit carbon at maximum rate allowed (which is the current emission level $c$), whereas to the left of the curve, the
optimal policy is to immediately reduce the emission rate and move vertically downward
toward the curve.
\end{comment}

{\color{black} \noindent Let us now consider concrete numerical values for an illustration. We focus on magnitudes for large firms, where $X_t$ is measured in units of MtCO$_2$, $\mu$ in MtCO$_2$/year and $\sigma$ in MtCO$_2/\sqrt{\text{year}}$ (it would, alternatively, be tCO$_2$ for individuals and GtCO$_2$ for countries). For a maximal excess-emission rate of $\overline{c}$ MtCO$_2$/year, an initial budget $X_0=x\approx H\cdot \overline{c}$ refers to an intended depletion horizon of $H$ years if the entity keeps emitting at $\overline{c}$ (ignoring drift and volatility of $X_t$). Avoidable excess emissions of several MtCO$_2$/year may be considered reasonable (we choose $\overline{c}=2$ in the examples in the sequel), and time horizons $H$ may vary from a few to 25 years, resulting in respective values for the initial budget $X_0=x$ (note that $x$ is a variable in our analysis).\footnote{{\color{black}According to the International Energy Agency \cite{IEA}, the global energy-related CO$_2$ emissions were about 37.8 GtCO$_2$ in 2024, and Lamboll et al.\ \cite{Lamboll} quantified the remaining carbon budget for a 50\% chance of staying within the 1.5$^\circ$C scenario to around 250 GtCO$_2$ and within the 2$^\circ$C scenario to around 1,200 GtCO$_2$, respectively, providing a framework for choices of $H$ depending on ambitions and the considered scenario. For instance, $x=10$ for $\overline{c}=2$ would then roughly refer to a time horizon of $H=5$ years, somewhat in line with a proportional appropriate share for a global 1.5$^\circ$C target in the above sense. }} The drift $\mu$ reflects potential regulatory tightening, technological progress, revised allocation, offsets outside the model, demand growth, or reclassification of unavoidable emissions. An annual deterministic change of $\mu\in(-0.5,1)$ MtCO$_2$/year looks reasonable depending on the scenario under consideration. The volatility parameter $\sigma$ represents uncertainty in the effective excess budget, such as demand shocks, production shocks, regulatory shocks, permit-allocation uncertainty, measurement revisions, or technological uncertainty. Values in the interval $\sigma\in(0.5,3)$ MtCO$_2/\sqrt{\text{year}}$ seem realistic and we choose $\sigma=1$ in the sequel. For the choice of $\gamma$, choices in the range of 100 $\text{\euro}$/tCO$_2$ could be meaningful (assuming $c_{net}=200\text{\euro}$/tCO$_2$ and $c_{tax}=100\text{\euro}$/tCO$_2$). However, as mentioned in Section \ref{Model and basic results}, we will in any case, and without loss of generality, scale the problem to $\gamma=1$ by suitably adapting the sustainability parameter $\Lambda$. \\

Consider therefore the example $\mu=1$, $\sigma=1$, $q=0.1,\Lambda=1.5$ and $S=[0,2]$ (that is, $\overline{c}=2$). Figure \ref{Ex1_Str} shows $V^{S}(x,2)$ (solid line), as a function of the initial carbon budget $x$. The dashed curve in the figure represents the classical optimal value function $V_{D}(x,c)=V_{D}(x)$ of the unconstrained case (that is, in the absence of downward ratcheting, cf.\ Remark \ref{Optima sin ratcheting}) for these parameters. One can see that the “cost” imposed by the downward ratcheting constraint is relatively limited. Moreover, a policy of continuous abatement is psychologically easier to implement than a strategy of remaining fully greedy and then abruptly reducing excess emissions to zero whenever the carbon budget falls below the  fixed barrier—an approach that would maximize the value function in the unconstrained case. Figure \ref{Ex1_Val} depicts the optimal
abatement threshold as a function of current available excess carbon emission level $x$. Since $\Lambda=1.5>\sqrt{\mu^{2}+2q\sigma^{2}}\approx 
1.095$, the optimal threshold is positive for all 
$0\leq c\leq2$  (cf.\ Proposition
\ref{Barrera cero}). % In other words, under the given model set-up, only for a value  $\Lambda>1.095$ one will reduce the emission rate xxx
Note that the optimal constant barrier in the absence of a ratcheting-down constraint is 2.997 for the present example (which one can for instance calculate with the formulas from Albrecher and Zhu \cite{AZ25}), plotted as a dashed line in Figure \ref{Ex1_Val}.\\

% \begin{example}\normalfont We consider in this example $\mu=3$, $\sigma=2$ and $q=0.1,$ $\Lambda=4$ for
% $S=[0,4]$. Figure \ref{Ex1_Str} shows $V^{S}(x,4)$ (solid line), as a function of
% initial carbon budget $x$. The dashed line in the figure represents the classical optimal value function $V_{\text{class}}(x,c)=V_{\text{class}}(x)$ of the unconstrained case (that is, in the absence of downward ratcheting) for these parameters. One can see that the “cost” imposed by the downward ratcheting constraint is relatively limited. Moreover, a policy of continuous abatement is psychologically easier to implement than a strategy of remaining fully greedy and then abruptly reducing excess emissions to zero whenever the carbon budget falls below the  fixed barrier—an approach that would maximize the value function in the unconstrained case. Figure \ref{Ex1_Val} depicts the optimal
% abatement threshold as a function of current available excess carbon emission level $x$. Since $\Lambda=4>\sqrt{\mu^{2}+2q\sigma^{2}}\approx 
% 3.1305$, the optimal threshold is positive for all 
% $0\leq c\leq4$  (cf.\ Proposition
% \ref{Barrera cero}). Note that the optimal constant barrier in the absence of a ratcheting-down constraint is 4.682 for the present example (which one can for instance calculate with the formulas from \cite{AZ25}). \hfill $\diamond$

\begin{figure}[ptb]
	\begin{subfigure}{.45\textwidth}
		\centering
		\includegraphics[width=5.5cm]{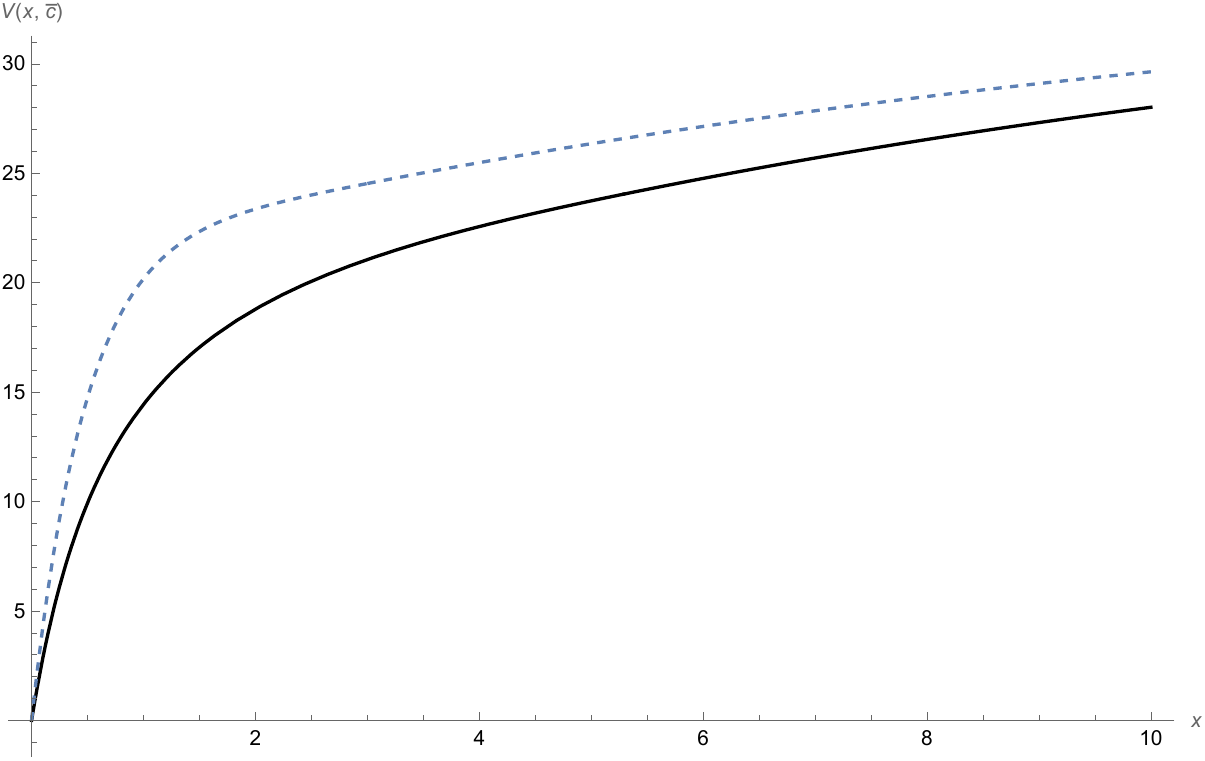}
		\caption{$V^S(x,2)$ and $V_{D}(x)$}
		\label{Ex1_Str}
	\end{subfigure} \begin{subfigure}{.45\textwidth}
		\centering
		\includegraphics[width=5.5cm]{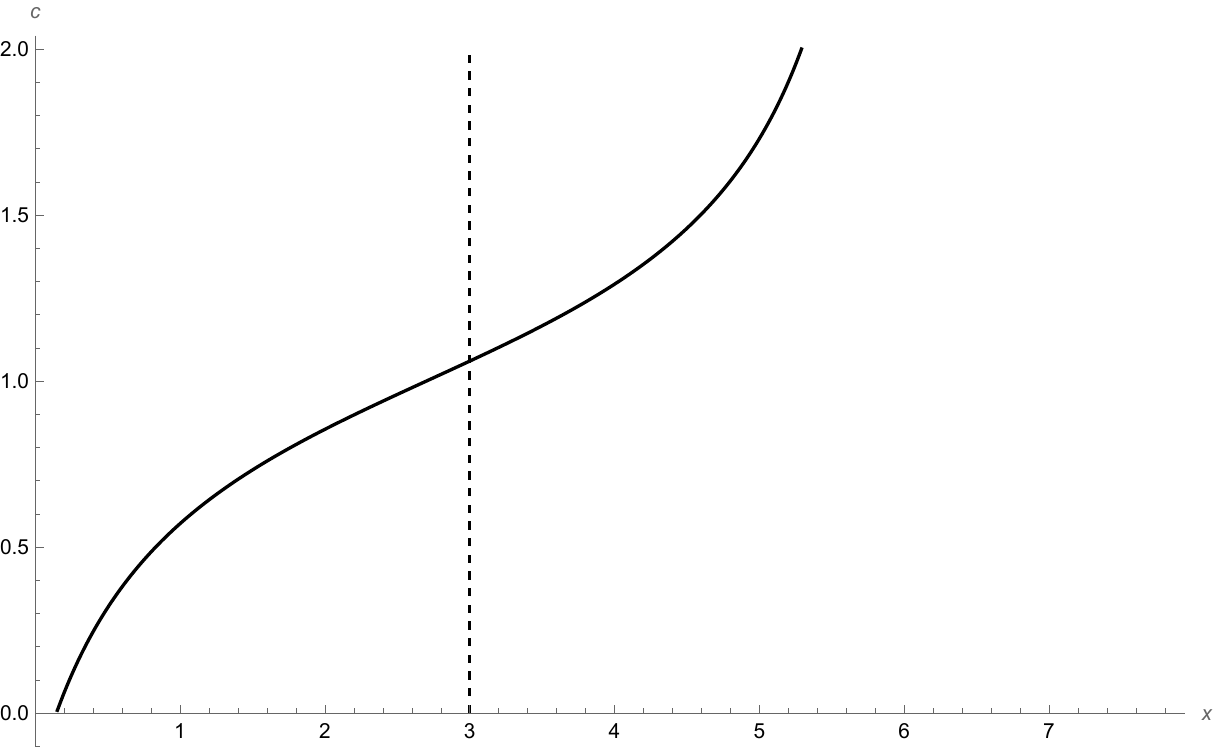}
		\caption{$z^*(c)$}
		\label{Ex1_Val}
	\end{subfigure}
	\caption{\protect\small Optimal value function $V^S(x,2)$ (solid line) and unconstrained value function $V_{D}(x)$ (dashed line) as well as optimal threshold $z^*(c)$ (right) for $\mu=1$, $\sigma=1$, $q=0.1$, $\Lambda=1.5$ and $S=[0,2]$.}%
	\label{fig1}%
\end{figure}

%\[%
%\begin{array}
%[c]{cc}%
%%TCIMACRO{\FRAME{itbpFU}{3.3693in}{2.1724in}{-0.0069in}{\Qcb{Fig ?.1 Ex 1}}%
%%{}{ex1optimalvnew.tiff}{\special{ language "Scientific Word";
%%type "GRAPHIC";  maintain-aspect-ratio TRUE;  display "USEDEF";
%%valid_file "F";  width 3.3693in;  height 2.1724in;  depth -0.0069in;
%%original-width 7.4996in;  original-height 4.8196in;  cropleft "0";
%%croptop "1";  cropright "1";  cropbottom "0";
%%filename 'Ex1OptimalVNew.tiff';file-properties "XNPEU";}} }%
%%BeginExpansion
%\raisebox{0.0069in}{\parbox[b]{3.3693in}{\begin{center}
%\includegraphics[
%natheight=4.819600in,
%natwidth=7.499600in,
%height=2.1724in,
%width=3.3693in
%]%
%{Ex1OptimalVNew.tiff}%
%\\
%Fig ?.1 Ex 1
%\end{center}}}
%%EndExpansion
%&
%%TCIMACRO{\FRAME{itbpFU}{2.9922in}{1.9623in}{0in}{\Qcb{Fig ?.2 Ex1 }}%
%%{}{ex1fronteranewsintecho.tiff}{\special{ language "Scientific Word";
%%type "GRAPHIC";  maintain-aspect-ratio TRUE;  display "USEDEF";
%%valid_file "F";  width 2.9922in;  height 1.9623in;  depth 0in;
%%original-width 7.4996in;  original-height 4.9026in;  cropleft "0";
%%croptop "1";  cropright "1";  cropbottom "0";
%%filename 'Ex1FronteraNewSinTecho.tiff';file-properties "XNPEU";}} }%
%%BeginExpansion
%{\parbox[b]{2.9922in}{\begin{center}
%\includegraphics[
%natheight=4.902600in,
%natwidth=7.499600in,
%height=1.9623in,
%width=2.9922in
%]%
%{Ex1FronteraNewSinTecho.tiff}%
%\\
%Fig ?.2 Ex1
%\end{center}}}
%%EndExpansion
%\end{array}
%\]

%\end{example}

\textbf{Sensitivity with respect to the drift $\mu$.} We now would like to focus on the effect of the drift $\mu$. 
%\begin{example}\normalfont We now would like to focus on the effect of the drift $\mu$ and choose $\sigma=1$, $q=0.1,\Lambda=1.5$ and $S=[0,2]$. 
}
Figure \ref{fig2} depicts $V^S(x,2)$ and the optimal abatement curve $z^*(c)$ for various values of $\mu$ (positive, zero and negative).
 In this example $\Lambda=1.5>\sqrt{\mu^{2}+2q\sigma^{2}}$ and $\Lambda+\mu>0$ for all chosen values of $\mu$. Therefore, by Remark \ref{rem33} and monotonicity of $z^*(c)$, the optimal threshold remains positive for all values of $c$. One observes empirically that the threshold function $z^*(c)$ is a convex function of $c$ for $c<c_e$ for some critical value $c_e$, and a concave function for $c>c_e$ (note that the $c$-axis is the ordinate in that plot). Moreover, the plot suggests that $c_e=\mu$ for any $\mu\ge 0$. We can not prove this latter claim with the techniques developed in this paper, but believe it to hold in general, and leave it as a conjecture for future research.   
 
 The optimal barrier without the abatement constraint is 2.997 for $\mu=1$, it is $4.059$ for $\mu=0.5$, 5.584 for $\mu=0$ and 6.110 for $\mu=-0.5$, respectively. Especially for zero or negative drift, this means that in the unconstrained case one would not allow carbon emissions unless the carbon budget level is quite high, as the budget would be depleted too quickly and the $\Lambda$-reward for the budget to last longer outweighs the immediate consumption benefit. Especially in such a situation, the abatement schedule is clearly preferable as it starts with consumption immediately and the efficiency loss (in terms of value function when compared to the unconstrained case) is still quite limited: for instance, for $x=5$ and $\mu=0$ the threshold strategy with optimal barrier 5.584 (emissions at rate $\overline{c}=2$ above the barrier and no emissions when the surplus is below the barrier) leads to a value function of 14.22, and for $x=5$ and $\mu=-0.5$ the threshold strategy with optimal barrier 6.110 leads to a value function of 8.71. In view of the values for $x=5$ in Figure \ref{Ex2_Str}, the efficiency loss from the non-ratcheting constraint is indeed quite small.    \\
 
 It is also a natural question to see by how much the optimal abatement strategy outperforms a simple (intuitive, but non-optimal)  linear abatement schedule $c(t)=\overline{c}-m\,t$ over time, starting in $c=\overline{c}$ and decreasing at a slope $m$ such that the original budget $x$ is used up when $c(t)$ hits $c=0$ (which we denote $t^*$; in case of $\sigma=0$ this would exactly mark the depletion time $\tau$).\footnote{{\color{black}Note that this is different from the schedule $c(t)=\overline{c}\, t$ without emission reduction that led to depletion horizon $H$ rather than $t^*$ before.}} A simple calculation gives $m= \overline{c}^2/(2 x)$ and $t^*=2x/\overline{c}$. {\color{black}This gives, due to the scheduled deterministic reduction, another simple benchmark relation between $\overline{c}$, initial budget $x$ and the envisaged time horizon without control (in this case $t^*$ instead of $H$). For instance, for an envisaged time horizon of $t^*=5$ years and $\overline{c}=2$, $x=5$ is now sufficient (to be compared with $x=10$ for $H=5$ years from before).} 
 %Since the initial budget $x$ will typically be given, a target time horizon for net-zero may determine the choice of the initial carbon emission rate $\overline{c}$. 
 %For $x=5$, the above choice $\overline{c}=2$ gives for instance $t^*=5$ years. 
 For the present parameters with $\mu=0$, a Monte Carlo simulation shows that such a {\color{black}simple} linear abatement schedule would lead to a value function of $9.81(\pm 0.14)$, where here and throughout, the number in parentheses indicates the halfwidth of the asymptotic 95\% confidence interval of the simulation.
  %(obtained for instance through a  or using an exact formula for a resulting consumption problem with linear decreasing drift $\mu(t)=0-m\,t=-0.4\,t$ for $X_t^C$), 
 %and the probability to have a depletion time larger than t=5 is only 0.33!
 %decomposition: C-part in V amounts to 3.72, Lambda part amounts to 6.14..
  The value of 9.81 is about 30\% below the corresponding value contained in the dotted line in Figure \ref{Ex2_Str}. Figure \ref{fig2a} illustrates both strategies for a sample path of the original surplus process $X_t$ for $\mu=0$. The black curve represents the surplus process $X_t^C$ when applying the optimal abatement strategy $z^*(c)$, which is the dotted line in Figure \ref{Ex2_Val}, and the black curve is the resulting abatement schedule for this sample path as a function of time. For this particular sample path, around $t=10$ the controlled carbon emission budget undershoots $x\approx 2.8$ for the first time, which leads to $c=0$ from then on, and the carbon budget remains positive for much longer (rewarded by the $\Lambda$-term, capturing the value of retaining part of the carbon budget, for example for the next generation). The red dotted curve is the (non-adaptive) linear abatement schedule $c(t)=2-0.4 \,t$, which equals zero (and therefore stopping emissions completely) already after 5 years, and the associated surplus process $X_t^C$ is depleted already much earlier than the one for the optimal strategy (a Monte Carlo simulation indicates that the expected depletion times under the two strategies are $38.50(\pm 1.37)$ and $9.83(\pm 0.70)$, respectively).\\ %\hfill $\diamond$ 

\begin{figure}[ptb]
	\begin{subfigure}{.45\textwidth}
		\centering
		\includegraphics[width=5.5cm]{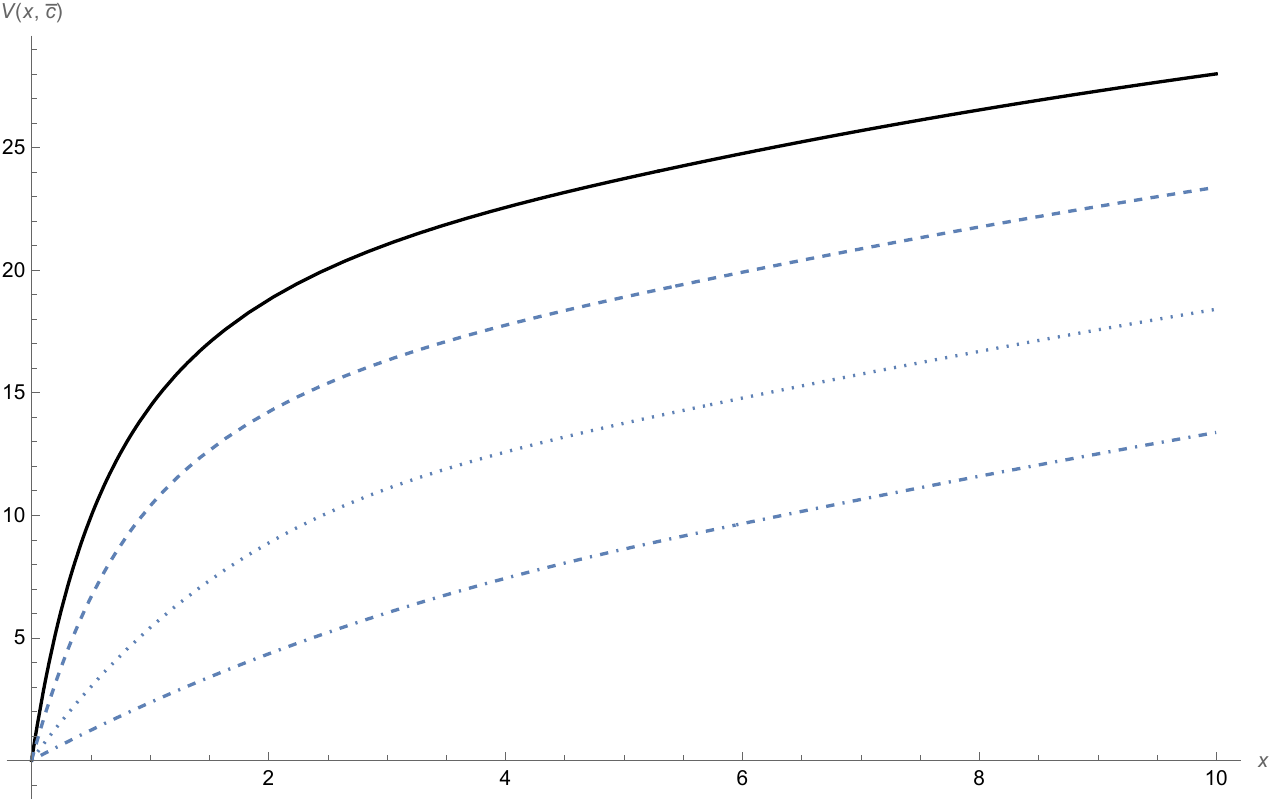}
		\caption{$V^S(x,2)$}
		\label{Ex2_Str}
	\end{subfigure} \begin{subfigure}{.45\textwidth}
		\centering
		\includegraphics[width=5.5cm]{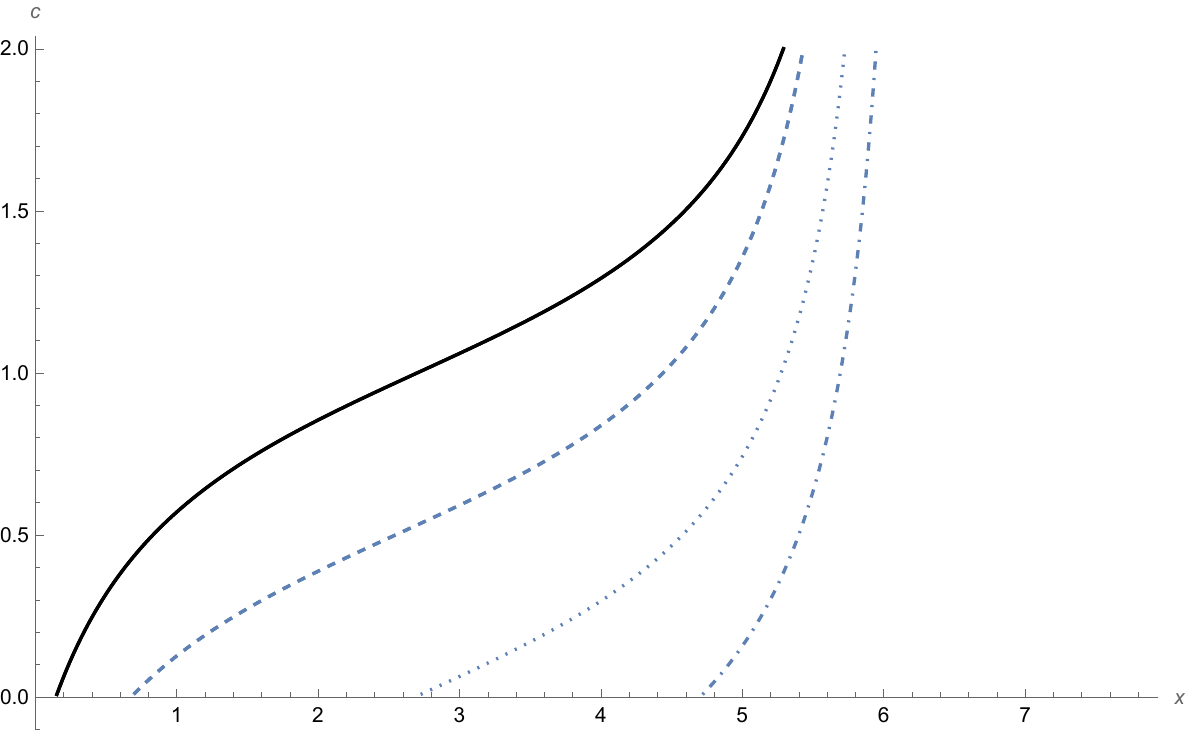}
		\caption{$z^*(c)$}
		\label{Ex2_Val}
	\end{subfigure}
	\caption{\protect\small Optimal value function $V^S(x,2)$ and optimal threshold $z^*(c)$ for $\sigma=1$, $q=0.1$, $\Lambda=1.5$ and $S=[0,2]$ for $\mu=1$ (solid line), $\mu=0.5$ (dashed line),  $\mu=0$ (dotted line) and $\mu=-0.5$ (dash-dotted line).
	}%
	\label{fig2}%
\end{figure}

\begin{figure}[ptb]
		\centering
		\includegraphics[width=10cm]{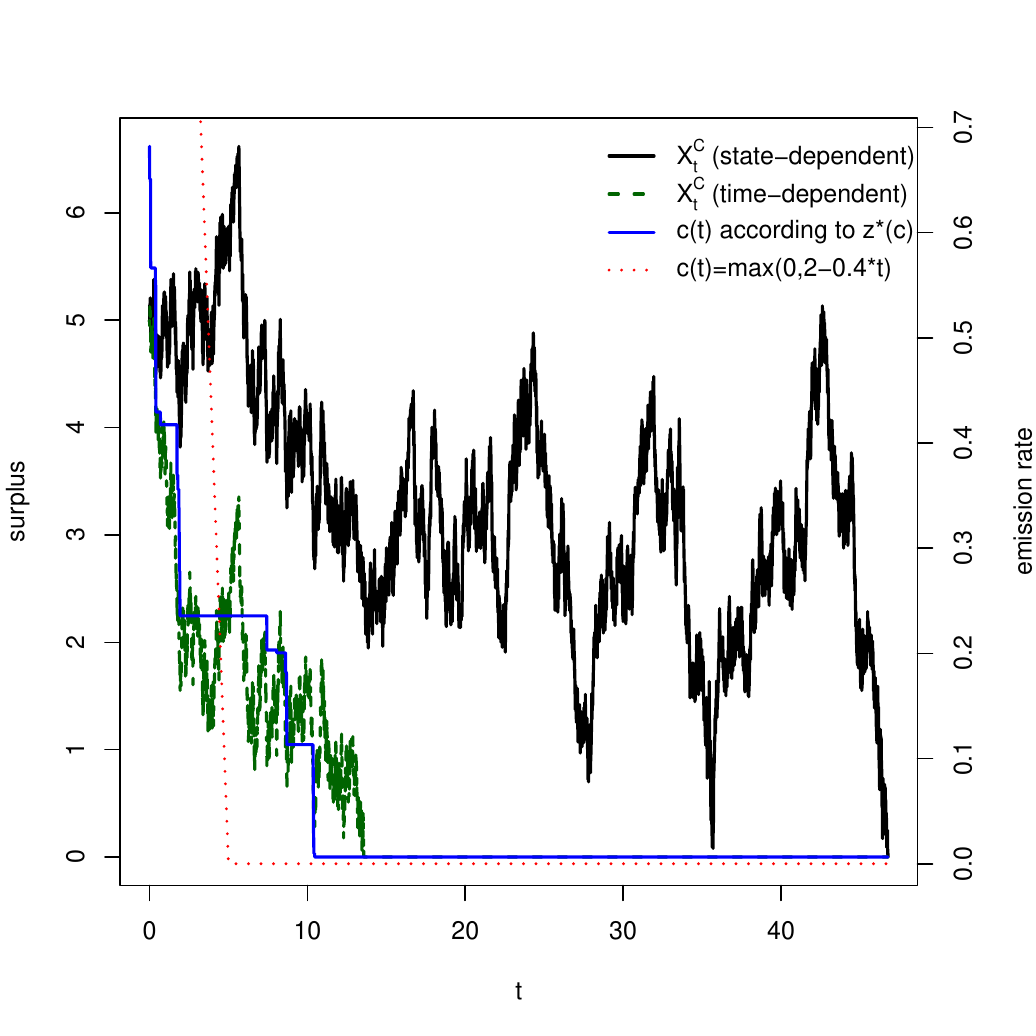}
		\caption{Sample path $X_t^C$ and resulting emission patterns for the parameters of Figure \ref{Ex2_Str} with $\mu=0$ for the optimal strategy according to $z^*(c)$ and a linear decreasing emission rate $c(t)=2-0.4t$.}
	\label{fig2a}%
\end{figure}
%\[%
%%TCIMACRO{\FRAME{itbpFU}{3.5345in}{2.2814in}{0in}{\Qcb{Fig ?.1 Ex 2}}%
%%{}{comparacionvaluefunctions.tiff}{\special{ language "Scientific Word";
%%type "GRAPHIC";  maintain-aspect-ratio TRUE;  display "USEDEF";
%%valid_file "F";  width 3.5345in;  height 2.2814in;  depth 0in;
%%original-width 7.4996in;  original-height 4.8196in;  cropleft "0";
%%croptop "1";  cropright "1";  cropbottom "0";
%%filename 'ComparacionValueFunctions.tiff';file-properties "XNPEU";}} }%
%%BeginExpansion
%{\parbox[b]{3.5345in}{\begin{center}
%\includegraphics[
%natheight=4.819600in,
%natwidth=7.499600in,
%height=2.2814in,
%width=3.5345in
%]%
%{ComparacionValueFunctions.tiff}%
%\\
%Fig ?.1 Ex 2
%\end{center}}}
%%EndExpansion
%\]%
%\[%
%%TCIMACRO{\FRAME{itbpFU}{3.9989in}{2.5201in}{0in}{\Qcb{Fig ?.2 Example 2}}%
%%{}{comparacionfronterasex2sintecho.tiff}%
%%{\special{ language "Scientific Word";  type "GRAPHIC";
%%maintain-aspect-ratio TRUE;  display "USEDEF";  valid_file "F";
%%width 3.9989in;  height 2.5201in;  depth 0in;  original-width 7.4996in;
%%original-height 4.708in;  cropleft "0";  croptop "1";  cropright "1";
%%cropbottom "0";
%%filename 'ComparacionFronterasEx2SinTecho.tiff';file-properties "XNPEU";}} }%
%%BeginExpansion
%{\parbox[b]{3.9989in}{\begin{center}
%\includegraphics[
%natheight=4.708000in,
%natwidth=7.499600in,
%height=2.5201in,
%width=3.9989in
%]%
%{ComparacionFronterasEx2SinTecho.tiff}%
%\\
%Fig ?.2 Example 2
%\end{center}}}
%%EndExpansion
%\]
%\end{example}
Note that the numerical value of $\Lambda$ balances the importance of substantial early emissions against the desire to delay the depletion time. It is therefore of interest to examine the sensitivity w.r.t. $\Lambda$ in mode detail. \\
%In the next example we therefore look at the sensitivity of the optimal strategy with respect to $\Lambda$. 

\begin{figure}[ptb]
	\begin{subfigure}{.45\textwidth}
		\centering
		\includegraphics[width=5.5cm]{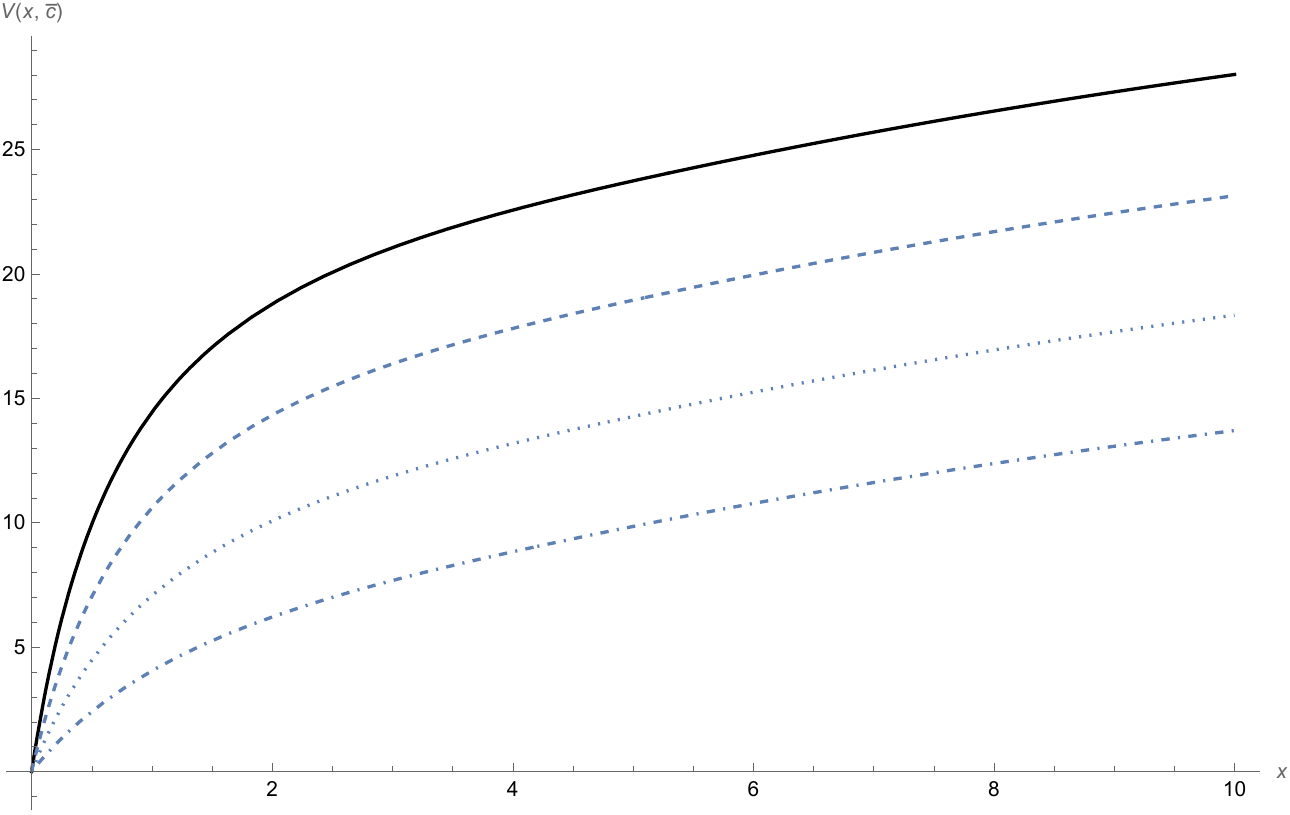}
		\caption{$V^S(x,2)$}
		\label{Ex3_Str}
	\end{subfigure} \begin{subfigure}{.45\textwidth}
		\centering
		\includegraphics[width=5.5cm]{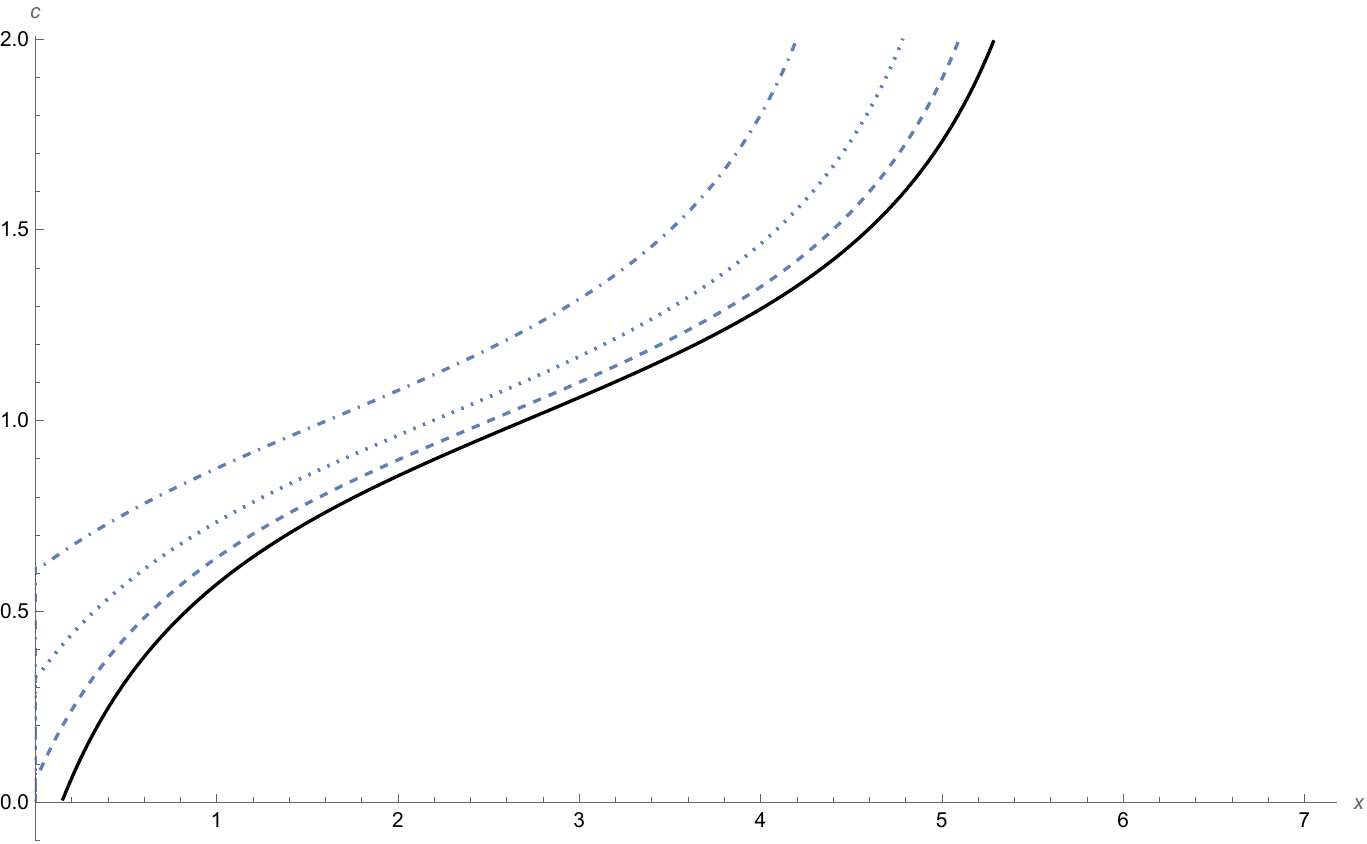}
		\caption{$z^*(c)$}
		\label{Ex3_Val}
	\end{subfigure}
	\caption{\protect\small Optimal value function $V^S(x,2)$ and optimal threshold $z^*$ as a function of $c$ for $\sigma=1$, $q=0.1$, $\mu=1$ and $S=[0,2]$ for $\Lambda=1.5$ (solid), $\Lambda=1$ (dashed), $\Lambda=0.5$ (dotted) and $\Lambda=0$ (dash-dotted).}%
	\label{fig3}%
\end{figure}

%\begin{example}\label{ex7.3}\normalfont 
{\color{black}\textbf{Sensitivity with respect to the reward parameter $\Lambda$.}}
%We now vary the reward level $\Lambda.$
For the case $\sigma=1$, $q=0.1,\mu=1$ and $S=[0,2]$, Figure \ref{fig3} depicts the value function and optimal threshold strategy for $\Lambda=1.5$ (solid line), $\Lambda=1$ (dashed line), $\Lambda=0.5$ (dotted line) and $\Lambda=0$ (dash-dotted line). For  $\Lambda=1.5>\sqrt{\mu^{2}+2q\sigma^{2}}\approx 
1.095$, the optimal strategy for each value of $c$ involves a
positive threshold. In contrast, for the smaller rewards $\Lambda=1,0.5$ and 0, respectively, the optimal threshold
is zero for all $0\leq c\leq\frac{2q\sigma^{2}+\mu^{2}-\Lambda^{2}}{2\left(
\Lambda+\mu\right)  }$ and positive otherwise (cf.\ Proposition \ref{Barrera cero} and Remark \ref{rem33}). 

As $\Lambda$ decreases, the importance of prolonged carbon-budget availability diminishes, leading to lower surplus threshold levels for emission rate reductions. For the extreme case $\Lambda=0$, corresponding to the absence of sustainability considerations, the limiting value of $c$ for which the threshold is positive becomes $(\mu^2+2q\sigma^2)/(2\mu)=0.6$ for the present parameters (cf.\  Remark \ref{rem33}), which is precisely the value at which $z^*(c)$ intersects the $c$-axis in Figure \ref{Ex3_Val}.

In general, one observes from Figure \ref{fig3} that the overall shape of the optimal curve $z^*(c)$ is relatively robust w.r.t.\ the choice of $\Lambda$ (Figure \ref{Ex3_Val}), but the value function itself is quite sensitive (Figure \ref{Ex3_Str}). In other words, for these parameter values the reward procedure (and thus the sustainability component) contributes substantially to the overall performance of the optimal strategy; however, the strategies themselves are relatively insensitive to increases or decreases in the reward procedure. %{\color{black} This, in particular, provides the magnitudes for which the purely profit-driven continued emission xxx}
\hfill $\diamond$

%\end{example}

\section{Conclusion}
In this paper, we solved the control problem of identifying the optimal abatement schedule for excess carbon emissions under a diffusion-type carbon budget, where the objective function consists of the expected discounted cumulative emissions together with a reward accrued as long as the carbon budget remains undepleted. We then implemented the proposed numerical procedure to compute the abatement schedule across several concrete examples and compared the results with both the unconstrained solution (i.e., without abatement) and a benchmark policy featuring a simple linear reduction in the consumption rate over time.  The results indicate that an optimal policy of gradual reduction in excess consumption entails only a moderate loss in the value function relative to the fully optimal emission schedule, which typically exhibits substantial fluctuations in the emission rate. These findings may inform the design of reduction pathways toward envisaged net-zero targets over fixed time horizons that are easier to implement from both psychological and practical perspectives. The numerical illustrations further reveal that the choice of the reward parameter $\Lambda$ has a significant impact on the resulting value function, while the optimal abatement schedule itself remains relatively robust. An interesting direction for future research is to refine the specification of the objective function for particular applications and to re-examine the associated optimal control problem. {\color{black}Also, it could be interesting to consider a variant of the problem, where the excess emission budget is not governed by an arithmetic Brownian motion, but by a jump process, where, e.g., negative jumps occur according to the arrival of stricter regulations and positive jumps occur according to innovations in carbon capture technology and sequestration techniques. Further variants of interest could include a stochastic mechanism for the dynamics of the tax rate $c_{tax}$ over time, which would introduce another random factor in the objective function itself. Also, an explicit modeling of positions in carbon certificates together with their random price dynamics in the present model setup would be an interesting future research direction. \\

Concerning the model specifications used in this paper, in the numerical illustrations, we observed} that the optimal threshold function $z^*(c)$ exhibits an inflection point which appears to lie exactly on the line $c=\mu$ for any drift parameter $\mu \ge 0$. We believe this to hold in general and pose it as a conjecture for future research. Furthermore, it can be interesting to see how a relaxation of the ratcheting-down constraint to a drawdown constraint (under which one might still increase the emission rate by a certain percentage of its current value) would influence the results. In the context of dividend optimization, drawdown constraints have been studied as a generalization of ratcheting-up restrictions (cf.\ Albrecher et al.\ \cite{AAM22}), and the resulting analysis proved to be highly non-trivial. We therefore expect that the corresponding analysis in the present setting will be very intricate as well.

%
%\[%
%%TCIMACRO{\FRAME{itbpFU}{3.9729in}{2.5624in}{0in}{\Qcb{Fig. ?.1 Ex 3}}%
%%{}{optimalvaluefunctionsej3_{e}j3b_{e}j2.tiff}%
%%{\special{ language "Scientific Word";  type "GRAPHIC";
%%maintain-aspect-ratio TRUE;  display "USEDEF";  valid_file "F";
%%width 3.9729in;  height 2.5624in;  depth 0in;  original-width 7.4996in;
%%original-height 4.8196in;  cropleft "0";  croptop "1";  cropright "1";
%%cropbottom "0";
%%filename 'OptimalValuefunctionsEj3_Ej3b_Ej2.tiff';file-properties "XNPEU";}}
%%}%
%%BeginExpansion
%{\parbox[b]{3.9729in}{\begin{center}
%\includegraphics[
%natheight=4.819600in,
%natwidth=7.499600in,
%height=2.5624in,
%width=3.9729in
%]%
%{OptimalValuefunctionsEj3_Ej3b_Ej2.tiff}%
%\\
%Fig. ?.1 Ex 3
%\end{center}}}
%%EndExpansion
%\]
%%
%
%\[%
%%TCIMACRO{\FRAME{itbpFU}{3.9773in}{2.4984in}{0in}{\Qcb{Fig ?.2 Ex. 3}}%
%%{}{fronterasex3sintecho.tiff}{\special{ language "Scientific Word";
%%type "GRAPHIC";  maintain-aspect-ratio TRUE;  display "USEDEF";
%%valid_file "F";  width 3.9773in;  height 2.4984in;  depth 0in;
%%original-width 7.4996in;  original-height 4.6942in;  cropleft "0";
%%croptop "1";  cropright "1";  cropbottom "0";
%%filename 'fronterasEx3SinTecho.tiff';file-properties "XNPEU";}} }%
%%BeginExpansion
%{\parbox[b]{3.9773in}{\begin{center}
%\includegraphics[
%natheight=4.694200in,
%natwidth=7.499600in,
%height=2.4984in,
%width=3.9773in
%]%
%{fronterasEx3SinTecho.tiff}%
%\\
%Fig ?.2 Ex. 3
%\end{center}}}
%%EndExpansion
%\]
%

\section{Appendix}

\noindent\textit{Proof of Proposition \ref{Proposicion Viscosidad}.} The proof is an adapted version of the proof of Albrecher et al.\ \cite[Prop.3.1]{AAM1}, tailored to the present situation of down-ratcheting (for self-containedness we give it in its complete form here again). Let us show first that $V$ is a viscosity supersolution in $(0,\infty
)\times\lbrack0,\overline{c})$. By Proposition
\ref{Monotone Optimal Value Function}, $\partial_{c}V(x,c)\geq0$ in $(0,\infty
)\times\lbrack0,\overline{c})$ in the viscosity sense.

Consider $(x,c)\in(0,\infty)\times\lbrack0,\overline{c}]$ and the
admissible strategy $C\in\Pi_{x,c}^{S}$, which emits at constant rate $c$ up
to the depletion time $\tau$. Let $X_{t}^{C}$ be the corresponding controlled
surplus process and suppose that there exists a test function $\varphi$ for
supersolution (\ref{HJB equation}) at $(x,c),$ then $\varphi\leq V$ and
$\varphi(x,c)=V(x,c)$. We want to prove that $\mathcal{L}^{c}\mathcal{(}%
\varphi)(x,c)\leq0$. For that purpose, we consider an auxiliary test function
for the supersolution $\tilde{\varphi}$ in such a way that $\tilde{\varphi
}\leq\varphi\leq V$ in $[0,\infty)\times\lbrack0,\overline{c}]$,
$\tilde{\varphi}=\varphi$ in $[0,2x]\times\lbrack0,\overline{c}]$ (so
$\mathcal{L}^{c}\mathcal{(}\varphi)(x,c)=\mathcal{L}^{c}\mathcal{(}%
\tilde{\varphi})(x,c)$) and $\mathcal{L}^{c}\mathcal{(}\tilde{\varphi}%
)(\cdot,c)$ is bounded in $[0,\infty)$. We introduce $\tilde{\varphi}$ because
$\mathcal{L}^{c}\mathcal{(}\varphi)(\cdot,c)$ could be unbounded in
$[0,\infty)$. We construct $\tilde{\varphi}$ as follows: take $g:[0,\infty
)\rightarrow\lbrack0,1]$ twice continuously differentiable with $g=0$ in
$[2x+1,\infty)$ and $g=1$ in $[0,2x]$, and define $\tilde{\varphi}(y,d)=$
$\varphi(y,d)g(y)$.
Using Lemma \ref{DPP}, we obtain for $h>0$%

\[%
\begin{array}
[c]{lll}%
\tilde{\varphi}(x,c) & = & V(x,c)\\
& \geq & \mathbb{E}\left[  \int\nolimits_{0}^{\tau\wedge h}e^{-q\,s}%
\,(c+\Lambda)ds\right]  +\mathbb{E}\left[  e^{-q(\tau\wedge h)}\tilde{\varphi
}(X_{\tau\wedge h}^{C},c)\right]  \text{.}%
\end{array}
\]
Hence, we get, using It\^{o}'s formula,%
\[%
\begin{array}
[c]{lll}%
0 & \geq & \mathbb{E}\left[  \int\nolimits_{0}^{\tau\wedge h}e^{-q\,s}%
\,(c+\Lambda)ds\right]  +\mathbb{E}\left[  e^{-q(\tau\wedge h)}\tilde{\varphi
}(X_{\tau\wedge h}^{C},c)-\tilde{\varphi}(x,c)\right] \\
& = & \mathbb{E}\left[  \int\nolimits_{0}^{\tau\wedge h}e^{-q\,s}%
\,(c+\Lambda)ds\right]  +\mathbb{E}\left[  \int\nolimits_{0}^{\tau\wedge
h}e^{-q\,s}(\frac{\sigma^{2}}{2}\partial_{xx}\tilde{\varphi}(X_{s}%
^{C},c)+\partial_{x}\tilde{\varphi}(X_{s}^{C},c)(\mu-c)-q\tilde{\varphi}%
(X_{s}^{C},c))ds\right] \\
&  & +\mathbb{E}\left[  \int_{0}^{\tau\wedge h}\partial_{x}\tilde{\varphi
}(X_{s}^{C},c)\sigma dWs~\right] \\
& = & \mathbb{E}\left[  \int\nolimits_{0}^{\tau\wedge h}e^{-q\,s}%
\mathcal{L}^{c}\mathcal{(}\tilde{\varphi})(X_{s}^{C},c)ds\right]  \text{.}%
\end{array}
\]
Since $\tau>0$ a.s.,%
\[
\left\vert \frac{1}{h}\int\nolimits_{0}^{\tau\wedge h}e^{-q\,s}\mathcal{L}%
^{c}\mathcal{(}\tilde{\varphi})(X_{s}^{C},c)ds\right\vert \leq\sup
_{y\in\lbrack0,\infty)}\left\vert \mathcal{L}^{c}\mathcal{(}\tilde{\varphi
})(y,c)\right\vert ,
\]
and
\[
\lim_{h\rightarrow0^{+}}\frac{1}{h}\int\nolimits_{0}^{\tau\wedge h}%
e^{-q\,s}\mathcal{L}^{c}\mathcal{(}\tilde{\varphi})(X_{s}^{C}%
,c)ds=\mathcal{L(}\tilde{\varphi})(x,c)~\text{a.s.;}%
\]
we conclude, using the bounded convergence theorem, that $\mathcal{L}%
^{c}\mathcal{(}\varphi)(x,c)=\mathcal{L}^{c}\mathcal{(}\tilde{\varphi
})(x,c)\leq0$; so $V$ is a viscosity supersolution at $(x,c)$.

%*****
%
%Difference with ratcheting:In the next proof we consider $[c-h,c]$ instead of
%$[c,c+h],C_{t}$ is non-increasing \ instead of non-decreasing
%
%***********************

Let us prove now that $V$ is a viscosity subsolution in $(0,\infty
)\times\lbrack 0,\overline{c})$. Assume first that $V$ is not a
subsolution of (\ref{HJB equation}) at $\left(  x,c\right)  \in(0,\infty
)\times(0,\overline{c}]$. Then there exist $\varepsilon>0$, $0<h<\min\left\{
x/2,c/2\right\}  $ and a (2,1)-differentiable function $\psi$ with
$\psi(x,c)=V(x,c)$ such that $\psi\geq V$,%

\begin{equation}
\max\{\mathcal{L}^{c}(\psi)(y,d),-\partial_{c}\psi(y,d)\}\leq-q\varepsilon<0
\label{Desig1aprima}%
\end{equation}
for $\left(  y,d\right)  \in$ $[x-h,x+h]\times\lbrack c-h,c]$ and%

\begin{equation}
V(y,d)\leq\psi(y,d)-\varepsilon\label{Desig2a}%
\end{equation}
for $\left(  y,d\right)  \notin\lbrack x-h,x+h]\times\lbrack c-h,c]$. Consider
the controlled risk process $X_{t}$ corresponding to an admissible strategy
$C\in\Pi_{x,c}^{S}$ and define
\[
\tau^{\ast}=\inf\{t>0:\text{ }\left(  X_{t},C_{t}\right)  \notin\lbrack
x-h,x+h]\times\lbrack c-h,c]\}\text{.}%
\]
Since $C_{t}$ is non-increasing and right-continuous, it can be written as
\begin{equation}
C_{t}=c+\int\nolimits_{0}^{t}dC_{s}^{co}-\sum_{\substack{C_{s}<C_{s^{-}%
}\\0\leq s\leq t}}(C_{s^{-}}-C_{s}), \label{Lt para ITO}%
\end{equation}
where $C_{s}^{co}$ is a continuous and non-increasing function.

Take a (2,1)-differentiable function $\psi:(0,\infty)\times\lbrack
 0,\overline{c}]\rightarrow\lbrack0,\infty)$. Note that, by the
mean value theorem, we have in the case $C_{s}<C_{s^{-}}$ that there
exists ${c}_{s}\in(C_{s},C_{s^{-}})$ with%
\[
\psi(X_{s}^{C},C_{s^{-}})-\psi(X_{s}^{C},C_{s})=(C_{s^{-}}-C_{s})\partial
_{c}\psi(X_{s}^{C},{c}_{s}).
\]
Using the expression (\ref{Lt para ITO}) and the change of variables formula
(see for instance Protter \cite{Protter}), we can write
\begin{equation}%
\begin{array}
[c]{l}%
e^{-q\tau^{\ast}}\psi(X_{\tau^{\ast}}^{C},C_{^{^{^{(\tau^{\ast})^{-}}}}}%
)-\psi(x,c)\\%
\begin{array}
[c]{ll}%
= & \int\nolimits_{0}^{\tau^{\ast}}e^{-qs}\partial_{x}\psi(X_{s}^{C},C_{s^{-}%
})(\mu-C_{s^{-}})ds+\int\nolimits_{0}^{\tau^{\ast}}e^{-qs}\partial_{c}%
\psi(X_{s}^{C},C_{s^{-}})dC_{s}^{co}\\
& -\sum_{\substack{C_{s}<C_{s^{-}}\\0\leq s<\tau^{\ast}}}e^{-qs}(C_{s^{-}%
}-C_{s})\partial_{c}\psi(X_{s}^{C},{c}_{s})\\
& +\int\nolimits_{0}^{\tau^{\ast}}e^{-qs}(-q\psi(X_{s}^{C},C_{s^{-}}%
)+\frac{\sigma^{2}}{2}\partial_{xx}\psi(X_{s}^{C},C_{s^{-}}))ds+\int
\nolimits_{0}^{\tau^{\ast}}e^{-qs}\partial_{x}\psi(X_{s}^{C},C_{s^{-}})\sigma
dW_{s}.
\end{array}
\end{array}
\label{Paso 1}%
\end{equation}
Hence, using (\ref{Desig1aprima}) and that $c_{s}\in\lbrack C_{s},C_{s^{-}%
}]\subset\lbrack c-h,c]$ for $s\in\lbrack0,\tau^{\ast})$, we can write%
\[%
\begin{array}
[c]{l}%
\mathbb{E}\left[  e^{-q\tau^{\ast}}\psi(X_{\tau^{\ast}}^{C},C_{^{^{^{(\tau
^{\ast})^{-}}}}})\right]  -\psi(x,c)\\%
\begin{array}
[c]{ll}%
= & \mathbb{E}\left[  \int\nolimits_{0}^{\tau^{\ast}}e^{-qs}\mathcal{L}%
^{C_{s^{-}}}(\psi)(X_{s}^{C},C_{s^{-}})ds-\int\nolimits_{0}^{\tau^{\ast}%
}e^{-qs}(C_{s^{-}}+\Lambda)ds\right] \\
& +\mathbb{E}\left[  \int\nolimits_{0}^{\tau^{\ast}}e^{-qs}\partial_{c}%
\psi(X_{s}^{C},C_{s^{-}})dC_{s}^{c}-\sum_{\substack{C_{s}\neq C_{s^{-}}\\0\leq
s<\tau^{\ast}}}e^{-qs}(C_{s^{-}}-C_{s})\partial_{c}\psi(X_{s}^{C}%
,c_{s})\right] \\
%\leq & -q\varepsilon\mathbb{E[}\int\nolimits_{0}^{\tau^{\ast}}e^{-qs}%
%]-\int\nolimits_{0}^{\tau^{\ast}}e^{-qs}(C_{s^{-}}+\Lambda)ds+\int
%\nolimits_{0}^{\tau^{\ast}}e^{-qs}q\varepsilon dC_{s}^{c}-\sum
%_{\substack{C_{s}\neq C_{s^{-}}\\0\leq s<\tau^{\ast}}}e^{-qs}(C_{s^{-}}%
%-C_{s})q\varepsilon(\ast)\\
\leq & \mathbb{E}\left[  \varepsilon\left(  e^{-q\tau^{\ast}}-1\right)
-\int\nolimits_{0}^{\tau^{\ast}}e^{-qs}(C_{s^{-}}+\Lambda)ds+q\varepsilon
\left(  \int\nolimits_{0}^{\tau^{\ast}}e^{-qs}dC_{s}\right)  \right]  .
\end{array}
\end{array}
\]

From (\ref{Desig2a}) and using that $V$ is a function that is non-decreasing in
the second variable as well as that $C_{s}$ is a non-increasing process,
\[%
\begin{array}
[c]{l}%
\mathbb{E}\left[  e^{-q\tau^{\ast}}V(X_{\tau^{\ast}}^{C},C_{\tau^{\ast}%
})\right] \\%
\begin{array}
[c]{cl}%
%= & \mathbb{E}\left[  e^{-q\tau^{\ast}}V(X_{\tau^{\ast}}^{C},C_{\tau^{\ast}%
%})\right]  -\mathbb{E}\left[  e^{-q\tau^{\ast}}V(X_{\tau^{\ast}}%
%^{C},C_{^{^{^{(\tau^{\ast})^{-}}}}})\right]  +\mathbb{E}\left[  e^{-q\tau
%^{\ast}}V(X_{\tau^{\ast}}^{C},C_{^{^{^{(\tau^{\ast})^{-}}}}})\right]  (\ast)\\
%\leq & \mathbb{E}\left[  e^{-q\tau^{\ast}}\left(  V(X_{\tau^{\ast}}%
%^{C},C_{\tau^{\ast}})-V(X_{\tau^{\ast}}^{C},C_{^{^{^{(\tau^{\ast})^{-}}}}%
%}\right)  \right]  +\mathbb{E}\left[  e^{-q\tau^{\ast}}(\psi(X_{\tau^{\ast}%
%}^{C},C_{^{^{^{(\tau^{\ast})^{-}}}}})-\varepsilon)\right]  (\ast)\\
\leq & \mathbb{E}\left[  e^{-q\tau^{\ast}}\left(  V(X_{\tau^{\ast}}%
^{C},C_{\tau^{\ast}})-V(X_{\tau^{\ast}}^{C},C_{^{^{^{(\tau^{\ast})^{-}}}}%
}\right)  \right]  +\mathbb{E}\left[  \psi(x,c)-e^{-q\tau^{\ast}}%
\varepsilon\right]  +\mathbb{E}\left[  \psi(X_{\tau^{\ast}}^{C},C_{(\tau
^{\ast})^{-}})e^{-q\tau^{\ast}}-\psi(x,c)\right] \\
%\leq & 0+\mathbb{E}\left[  \psi(x,c)-e^{-q\tau^{\ast}}\varepsilon\right]
%+\mathbb{E}\left[  \varepsilon\left(  e^{-q\tau^{\ast}}-1\right)
%-\int\nolimits_{0}^{\tau^{\ast}}e^{-qs}(C_{s^{-}}+\Lambda)ds+q\varepsilon
%\left(  \int\nolimits_{0}^{\tau^{\ast}}e^{-qs}dC_{s}\right)  \right]  (\ast)\\
\leq & \psi(x,c)-\varepsilon-\mathbb{E}(\int\nolimits_{0}^{\tau^{\ast}}%
e^{-qs}(C_{s^{-}}+\Lambda)ds).
\end{array}
\end{array}
\]
Hence, using Lemma \ref{DPP}, we have that
\[
V(x,c)=\sup\limits_{C\in\Pi_{x,c}^{S}}\mathbb{E}\left(  \int\nolimits_{0}%
^{\tau^{\ast}}e^{-qs}(C_{s^{-}}+\Lambda)ds+e^{-q\tau^{\ast}}V(X_{\tau^{\ast}}%
^{C},C_{\tau^{\ast}})\right)  \leq\psi(x,c)-\varepsilon.
\]
But the latter is a contradiction because we have assumed that $V(x,c)=\psi(x,c)$.
When $c=0$, $V(x,0)$ solves $\mathcal{L}^{0}(V)(x,0)=0$, which gives the
result.\hfill $\blacksquare$\\

\noindent\textit{Proof of Lemma \ref{Lema para Unicidad}.}
 A locally Lipschitz function $\overline{u}$\ $:[0,\infty)\times\lbrack
0,\overline{c}]\rightarrow{\mathbb{R}}$\ is a viscosity supersolution of
(\ref{HJB equation}) at $(x,c)\in(0,\infty)\times(0,\overline{c})$, if any
test function $\varphi$ for supersolution at $(x,c)$ satisfies
\begin{equation}
\max\{\mathcal{L}^{c}(\varphi)(x,c),-\partial_{c}\varphi(x,c)\}\leq0\text{,}
\label{superyalfa}%
\end{equation}
and a locally Lipschitz function $\underline{u}:[0,\infty)\times
\lbrack0,\overline{c}]\rightarrow{\mathbb{R}}$\ is a viscosity subsolution\ of
(\ref{HJB equation}) at $(x,c)\in(0,\infty)\times(0,\overline{c})$\ if any
test function $\psi$ for subsolution at $(x,c)$ satisfies
\begin{equation}
\max\{\mathcal{L}^{c}(\psi)(x,c),-\partial_{c}\psi(x,c)\}\geq0.
\label{subxalfa}%
\end{equation}
Suppose that there is a point $(x_{0},c_{0})\in\lbrack0,\infty)\times
(0,\overline{c})$ such that $\underline{u}(x_{0},c_{0}%
)-\overline{u}(x_{0},c_{0})>0$. Let us define $h(c)=1+e^{{c}%
/{{\overline{c}}}}$ and
\[
\overline{u}^{s}(x,c)=s\,h(c)\,\overline{u}(x,c)
\]
for any $s>1$. We have that $\varphi$ is a test function for supersolution of
$\overline{u}$ at $(x,c)$ if and only if $\varphi^{s}=s\,h(c)\,\varphi$ is a
test function for supersolution of $\overline{u}^{s}$ at $(x,c)$. By
(\ref{superyalfa}) and using $1-s\,h(c)<1-s<0,$ we have {%
\begin{equation}%
\begin{array}
[c]{ccl}%
\mathcal{L}^{c}(\varphi^{s})(x,c) & = & \frac{\sigma^{2}}{2}\,s\,h(c)\,\partial
_{xx}\varphi\left(  x,c\right)  +(\mu-c)\,s\;h(c)\,\partial_{x}\varphi\left(
x,c\right)  -qs\,h(c)\varphi\left(  x,c\right)  +c+\Lambda\\
& = & s\,h(c)\,\mathcal{L}^{c}(\varphi)(x,c)+(c+\Lambda)(1-s\,h(c))\\
& < & 0
\end{array}
\label{Desigualdad L1}%
\end{equation}
} and%
\begin{equation}
\partial_{c}\varphi^{s}(x,c)\geq\frac{s}{\overline{c}}\,e^{{c}/{\overline
{c}}}\varphi(x,c)>0 \label{Desigualdad L2}%
\end{equation}
for $\varphi(x,c)>0$. Take $s_{0}>1$, then $\underline{u}(x_{0},c_{0}%
)-\overline{u}^{s}(x_{0},c_{0})>0$. We define%
\begin{equation}
M=\sup\limits_{x\geq0,0\leq c\leq\overline{c}}\left(  \underline
{u}(x,c)-\overline{u}^{s_{0}}(x,c)\right)  . \label{Definicion de M}%
\end{equation}
Since $\lim_{x\rightarrow\infty}\underline{u}(x,c)\leq(\overline{c}%
+\Lambda)/q\leq\lim_{x\rightarrow\infty}\overline{u}(x,c)$, there exists a 
$b>x_{0}$ such that%

\begin{equation}
\sup\limits_{0\leq c\leq\overline{c}}\underline{u}(x,c)-\overline{u}^{s_{0}%
}(x,c)<0\text{ for }x\geq b. \label{desborde}%
\end{equation}
From (\ref{desborde}), we obtain that%

\begin{equation}
0<\underline{u}(x_{0},c_{0})-\overline{u}^{s_{0}}(x_{0},c_{0})\leq
M:=\max\limits_{x\in\left[  0,b\right]  ,0\leq c\leq\overline{c}}\left(
\underline{u}(x,c)-\overline{u}^{s_{0}}(x,c)\right)  . \label{desmax}%
\end{equation}
Call $\left(  x^{\ast},c^{\ast}\right)  :=\arg\max\limits_{x\in\left[
0,b\right]  ,0\leq c\leq\overline{c}}\left(  \underline{u}(x,c)-\overline
{u}^{s_{0}}(x,c)\right)  $. Let us consider the set%

\[
\mathcal{A}=\left\{  \left(  x,y,c,d\right)  :0\leq x\leq y\leq b\text{,
}0\leq\ c\leq\overline{c}\text{, }0\leq d\leq\overline{c}\right\}
\]
and, for all $\lambda>0$, the functions%

\begin{equation}%
\begin{array}
[c]{l}%
\Phi^{\lambda}\left(  x,y,c,d\right)  =\dfrac{\lambda}{2}\left(  x-y\right)
^{2}+\dfrac{\lambda}{2}\left(  c-d\right)  ^{2}+\frac{2m}{\lambda^{2}\left(
y-x\right)  +\lambda},\\
\Sigma^{\lambda}\left(  x,y,c,d\right)  =\underline{u}(x,c)-\overline
{u}^{s_{0}}(y,d)-\Phi^{\lambda}\left(  x,y,c,d\right)  .
\end{array}
\label{def-sigma-definitiva-por-hoy}%
\end{equation}
$\allowbreak$ Calling $M^{\lambda}=\max\limits_{A}\Sigma^{\lambda}$ and
$\left(  x_{\lambda},y_{\lambda},c_{\lambda},d_{\lambda}\right)  =\arg
\max\limits_{A}\Sigma^{\lambda}$, we obtain that
\[
M^{\lambda}\geq\Sigma^{\lambda}(x^{\ast},x^{\ast},c^{\ast},c^{\ast}%
)=M-\frac{2m}{\lambda},%
\]
and so%

\begin{equation}
\liminf\limits_{\lambda\rightarrow\infty}M^{\lambda}\geq M.
\label{liminf-mlambda}%
\end{equation}
There exists $\lambda_{0}$ large enough and $s$ small enough such that if
$\lambda\geq\lambda_{0}$, then $\left(  x_{\lambda},y_{\lambda},c_{\lambda
},d_{\lambda}\right)  $ $\notin\partial A$, the proof is similar to the one of
Lemma 4.5 of Albrecher et al.\ {\cite{AAM}. }Using the inequality%

\[
\Sigma^{\lambda}\left(  x_{\lambda},x_{\lambda},c_{\lambda},c_{\lambda
}\right)  +\Sigma^{\lambda}\left(  y_{\lambda},y_{\lambda},d_{\lambda
},d_{\lambda}\right)  \leq2\Sigma^{\lambda}\left(  x_{\lambda},y_{\lambda
},c_{\lambda},d_{\lambda}\right)  ,
\]
we obtain that%
\[
\lambda\left\Vert (x_{\lambda}-y_{\lambda},c_{\lambda}-d_{\lambda})\right\Vert
_{2}^{2}\leq\underline{u}(x_{\lambda},c_{\lambda})-\underline{u}(y_{\lambda
},d_{\lambda})+\overline{u}^{s_{0}}(x_{\lambda},c_{\lambda})-\overline
{u}^{s_{0}}(y_{\lambda},d_{\lambda})+4m(y_{\lambda}-x_{\lambda}).
\]
Consequently
\begin{equation}
\lambda\left\Vert (x_{\lambda}-y_{\lambda},c_{\lambda}-d_{\lambda})\right\Vert
_{2}^{2}\leq6m\left\Vert (x_{\lambda}-y_{\lambda},c_{\lambda}-d_{\lambda
})\right\Vert _{2}. \label{desxalfa-yalfac}%
\end{equation}
We can find a sequence $\lambda_{n}\rightarrow\infty$ such that $\left(
x_{\lambda_{n}},y_{\lambda_{n}},c_{\lambda_{n}},d_{\lambda_{n}}\right)
\rightarrow\left(  \widehat{x},\widehat{y},\widehat{c},\widehat{d}\right)  \in
A$. From (\ref{desxalfa-yalfac}), we get that
\begin{equation}
\left\Vert (x_{\lambda_{n}}-y_{\lambda_{n}},c_{\lambda_{n}}-d_{\lambda_{n}
})\right\Vert _{2}\leq6m/\lambda_{n} , \label{cota con 6m}%
\end{equation}
which gives $\widehat{x}=\widehat{y}$ and $\widehat{c}=\widehat{d}$.

Since $\Sigma^{\lambda}\left(  x,y,c,d\right)  =\underline{u}(x,c)-\overline
{u}^{s_{0}}(y,d)-\Phi^{\lambda}\left(  x,y,c,d\right)  $ reaches the maximum
in $\left(  x_{\lambda},y_{\lambda},c_{\lambda},d_{\lambda}\right)  \ $in the
interior of the set $A,$ the function%
\[
\psi(x,c)=\Phi^{\lambda}\left(  x,y_{\lambda},c,d_{\lambda}\right)
-\Phi^{\lambda}\left(  x_{\lambda},y_{\lambda},c_{\lambda},d_{\lambda}\right)
+\underline{u}\left(  x_{\lambda},c_{\lambda}\right)
\]
is a test for subsolution for $\underline{u}$ of the HJB equation at the point
$\left(  x_{\lambda},c_{\lambda}\right)  $.

In addition, the function%
\[
\varphi^{s_{0}}(y,d)=-\Phi^{\lambda}\left(  x_{\lambda},y,c_{\lambda
},d\right)  +\Phi^{\lambda}\left(  x_{\lambda},y_{\lambda},c_{\lambda
},d_{\lambda}\right)  +\overline{u}^{s_{0}}\left(  y_{\lambda},d_{\lambda
}\right)
\]
is a test for supersolution for $\overline{u}^{s_{0}}$ at $\left(  y_{\lambda
},d_{\lambda}\right)  $ and so
\[
\partial_{c}\varphi^{s_{0}}(y_{\lambda},d_{\lambda})\geq\frac{s_{0}}{\overline{c}}\varphi(y_{\lambda},d_{\lambda}%
)e^{{c}/{\overline{c}}}>0,
\]
using $y_{\lambda}>0$. Consequently, $\partial_{c}\psi(x_{\lambda},c_{\lambda})=\partial_{c}\varphi^{s_{0}%
}(y_{\lambda},d_{\lambda})>0$, and so we have $\mathcal{L}^{c_{\lambda}}%
(\psi)(x_{\lambda},c_{\lambda})\geq0.$

\ Assume first that the functions $\underline{u}(x,c)$ and $\overline
{u}^{s_{0}}(y,d)$ are (2,1)-differentiable at $(x_{\lambda},c_{\lambda})\ $and
$(y_{\lambda},d_{\lambda})$ respectively. Since $\Sigma^{\lambda}$ defined in
(\ref{def-sigma-definitiva-por-hoy}) reaches a local maximum at $\left(
x_{\lambda},y_{\lambda},c_{\lambda},d_{\lambda}\right)  $ $\notin\partial A$,
we have that
\[
\partial_{x}\Sigma^{\lambda}\left(  x_{\lambda},y_{\lambda},c_{\lambda
},d_{\lambda}\right)  =\partial_{y}\Sigma^{\lambda}\left(  x_{\lambda
},y_{\lambda},c_{\lambda},d_{\lambda}\right)  =0
\]
and so%

\begin{equation}%
\begin{array}
[c]{lll}%
\partial_{x}\underline{u}(x_{\lambda},c_{\lambda}) & = & \partial_{x}
\Phi^{\lambda}(x_{\lambda},y_{\lambda},c_{\lambda},d_{\lambda})\\
& = & \lambda\left(  x_{\lambda}-y_{\lambda}\right)  +\frac{2m}{\left(
\lambda\left(  y_{\lambda}-x_{\lambda}\right)  +1\right)  ^{2}}\\
& = & -\partial_{y}\Phi^{\lambda}(x_{\lambda},y_{\lambda},c_{\lambda
},d_{\lambda})=\partial_{x}\overline{u}^{s_{0}}(y_{\lambda},d_{\lambda}).
\end{array}
\label{dos derivadas iguales}%
\end{equation}
Defining $A=\partial_{xx}\underline{u}(x_{\lambda},c_{\lambda})$ and
$B=\partial_{xx}\overline{u}^{s_{0}}(y_{\lambda},d_{\lambda})$, we obtain%

\[%
\begin{array}
[c]{l}%
\left(
\begin{array}
[c]{ll}%
\partial_{xx}\Sigma^{\lambda}\left(  x_{\lambda},y_{\lambda},c_{\lambda
},d_{\lambda}\right)  & \partial_{xy}\Sigma^{\lambda}\left(  x_{\lambda
},y_{\lambda},c_{\lambda},d_{\lambda}\right) \\
\partial_{xy}\Sigma^{\lambda}\left(  x_{\lambda},y_{\lambda},c_{\lambda
},d_{\lambda}\right)  & \partial_{yy}\Sigma^{\lambda}\left(  x_{\lambda
},y_{\lambda},c_{\lambda},d_{\lambda}\right)
\end{array}
\right) \\ \medskip
=\left(
\begin{array}
[c]{ll}%
A-\partial_{xx}\Phi^{\lambda}\left(  x_{\lambda},y_{\lambda},c_{\lambda
},d_{\lambda}\right)  & -\partial_{xy}\Phi^{\lambda}\left(  x_{\lambda
},y_{\lambda},c_{\lambda},d_{\lambda}\right) \\
-\partial_{xy}\Phi^{\lambda}\left(  x_{\lambda},y_{\lambda},c_{\lambda
},d_{\lambda}\right)  & -B-\partial_{yy}\Phi^{\lambda}\left(  x_{\lambda
},y_{\lambda},c_{\lambda},d_{\lambda}\right)
\end{array}
\right)  \leq0.
\end{array}
\]
It is hence a negative semi-definite matrix, and
\[%
\begin{pmatrix}
A & 0\\
0 & -B
\end{pmatrix}
\leq H\left(  \Phi^{\lambda}\right)  (x_{\lambda},y_{\lambda},c_{\lambda
},d_{\lambda}):=\left(
\begin{array}
[c]{ll}%
\partial_{xx}\Phi^{\lambda}\left(  x_{\lambda},y_{\lambda},c_{\lambda
},d_{\lambda}\right)  & \partial_{xy}\Phi^{\lambda}\left(  x_{\lambda
},y_{\lambda},c_{\lambda},d_{\lambda}\right) \\
\partial_{xy}\Phi^{\lambda}\left(  x_{\lambda},y_{\lambda},c_{\lambda
},d_{\lambda}\right)  & \partial_{yy}\Phi^{\lambda}\left(  x_{\lambda
},y_{\lambda},c_{\lambda},d_{\lambda}\right)
\end{array}
\right)  .
\]
%The latter matrix inequality entails that the difference is also a negative semi-definite matrix.

In the case that $\underline{u}(x,c)$ and $\overline{u}^{s_{0}}(y,d)$ are not
(2,1)-differentiable at $\left(  x_{\lambda},c_{\lambda}\right)  \ $and
$(y_{\lambda},d_{\lambda})$, respectively, we can resort to a more general
theorem to get a similar result. Using Theorem 3.2 of Crandall, Ishii and
Lions \cite{CrandallIshiLions}, it can be proved that for any $\delta>0$ 
there exist real numbers $A_{\delta}$ and $B_{\delta}$ such that
\begin{equation}%
\begin{pmatrix}
A_{\delta} & 0\\
0 & -B_{\delta}%
\end{pmatrix}
\leq H\left(  \Phi^{\lambda}\right)  (x_{\lambda},y_{\lambda},c_{\lambda
},d_{\lambda})+\delta\left(  H\left(  \Phi^{\lambda}\right)  (x_{\lambda
},y_{\lambda},c_{\lambda},d_{\lambda})\right)  ^{2} \label{MatrizDeltaruina}%
\end{equation}
and \bigskip%

\begin{equation}%
\begin{array}
[c]{c}%
\frac{\sigma^{2}}{2}A_{\delta}+(\mu-c_{\lambda})\partial_{x}\psi(x_{\lambda
},c_{\lambda})-q\psi(x_{\lambda},c_{\lambda})+c_{\lambda}+\Lambda\geq0,\\
\frac{\sigma^{2}}{2}B_{\delta}+(\mu-d_{\lambda})\partial_{x}\varphi^{s_{0}%
}(y_{\lambda},d_{\lambda})-q\varphi^{s_{0}}(y_{\lambda},d_{\lambda
})+d_{\lambda}+\Lambda\leq0.
\end{array}
\label{Inecuaciones con A y B delta}%
\end{equation}
The expression (\ref{MatrizDeltaruina}) implies that $A_{\delta}-B_{\delta
}\leq0$ because%
\[
H\left(  \Phi^{\lambda}\right)  (x_{\lambda},y_{\lambda},c_{\lambda
},d_{\lambda})=\partial_{xx}\Phi^{\lambda}\left(  x_{\lambda},y_{\lambda
},c_{\lambda},d_{\lambda}\right)
\begin{pmatrix}
1 & -1\\
-1 & 1
\end{pmatrix}
\]
and
\[
\left(  H\left(  \Phi^{\lambda}\right)  (x_{\lambda},y_{\lambda},c_{\lambda
},d_{\lambda})\right)  ^{2}=2\left(  \partial_{xx}\Phi^{\lambda}\left(
x_{\lambda},y_{\lambda},c_{\lambda},d_{\lambda}\right)  \right)  ^{2}%
\begin{pmatrix}
1 & -1\\
-1 & 1
\end{pmatrix}
.
\]
Therefore,%
\[%
\begin{array}
[c]{lll}%
A_{\delta}-B_{\delta} & = &
\begin{pmatrix}
1 & 1
\end{pmatrix}
\left(
\begin{array}
[c]{cc}%
A_{\delta} & 0\\
0 & -B_{\delta}%
\end{array}
\right)  \left(
\begin{array}
[c]{c}%
1\\
1
\end{array}
\right) \\
& \leq &
\begin{pmatrix}
1 & 1
\end{pmatrix}
\left(  H\left(  \Phi^{\lambda}\right)  (x_{\lambda},y_{\lambda},c_{\lambda
},d_{\lambda})+\delta\left(  H\left(  \Phi^{\lambda}\right)  (x_{\lambda
},y_{\lambda},c_{\lambda},d_{\lambda})\right)  ^{2}\right)  \left(
\begin{array}
[c]{c}%
1\\
1
\end{array}
\right) \\
& = & 0.
\end{array}
\]
And so, since $\varphi^{s_{0}}\left(  y_{\lambda},d_{\lambda}\right)
=\overline{u}^{s_{0}}\left(  y_{\lambda},d_{\lambda}\right)  $, $\psi
(x_{\lambda},c_{\lambda})=\underline{u}(x_{\lambda},c_{\lambda})$ and
\[
\partial_{x}\varphi^{s_{0}}\left(  y_{\lambda},d_{\lambda}\right)
=-\partial_{y}\Phi^{\lambda}\left(  x_{\lambda},y_{\lambda},c_{\lambda
},d_{\lambda}\right)  =\partial_{x}\Phi^{\lambda}\left(  x_{\lambda
},y_{\lambda},c_{\lambda},d_{\lambda}\right)  =\partial_{x}\psi(x_{\lambda
},c_{\lambda}),
\]
we obtain%

\begin{equation}%
\begin{array}
[c]{lll}%
\underline{u}(x_{\lambda},c_{\lambda})-\overline{u}^{s_{0}}\left(  y_{\lambda
},d_{\lambda}\right)  & = & \psi(x_{\lambda},c_{\lambda})-\varphi^{s_{0}
}\left(  y_{\lambda},d_{\lambda}\right) \\
& \leq & \frac{\sigma^{2}}{2q}(A_{\delta}-B_{\delta})\\
&  & +\left(  \frac{c_{\lambda}}{q}-\frac{d_{\lambda}}{q}\right)
(1-\partial_{x}\Phi^{\lambda}\left(  x_{\lambda},y_{\lambda},c_{\lambda
},d_{\lambda}\right)  )\\
& \leq & \left(  \frac{c_{\lambda}}{q}-\frac{d_{\lambda}}{q}\right)
\left(1-\lambda\left(  x_{\lambda}-y_{\lambda}\right)  -\frac{2m}{\left(
\lambda\left(  y_{\lambda}-x_{\lambda}\right)  +1\right)  ^{2}}\right).
\end{array}
\label{desdifsup1}%
\end{equation}
Hence, from (\ref{cota con 6m}) and (\ref{liminf-mlambda}), we get%

\begin{align*}
0  &  <M\leq\liminf\limits_{\lambda\rightarrow\infty}M_{\lambda}\leq
\lim\limits_{_{n\rightarrow\infty}}M_{\lambda_{n}}=\lim\limits_{_{n\rightarrow
\infty}}\Sigma^{\lambda_{n}}(x_{\lambda_{n}},y_{\lambda_{n}},c_{\lambda_{n}%
},d_{\lambda_{n}})=\underline{u}(\widehat{x},\widehat{c})-\overline{u}^{s_{0}%
}(\widehat{x},\widehat{c})\\
&  \leq\lim_{n\longrightarrow\infty}\left(  \frac{c_{\lambda_{n}}}{q}%
-\frac{d_{\lambda_{n}}}{q}\right)  (1-\lambda_{n}\left(  x_{\lambda_{n}%
}-y_{\lambda_{n}}\right)  -\frac{2m}{\left(  \lambda_{n}\left(  y_{\lambda
_{n}}-x_{\lambda_{n}}\right)  +1\right)  ^{2}})\\
&  \leq.\lim_{n\longrightarrow\infty}\left\vert \frac{c_{\lambda_{n}}}%
{q}-\frac{d_{\lambda_{n}}}{q}\right\vert (1+\lambda_{n}\left\Vert
(x_{\lambda_{n}}-y_{\lambda_{n}},c_{\lambda_{n}}-d_{\lambda_{n}})\right\Vert
_{2}+\frac{2m}{\left(  \lambda_{n}\left(  y_{\lambda_{n}}-x_{\lambda_{n}%
}\right)  +1\right)  ^{2}})\\
&  \leq\lim_{n\longrightarrow\infty}\left\vert \frac{c_{\lambda_{n}}}{q}%
-\frac{d_{\lambda_{n}}}{q}\right\vert (1+8m)=0.
\end{align*}
This is a contradiction and so we get the result.\hfill $\blacksquare$\\

\noindent\textit{Proof of Theorem \ref{Optima Discreta Threshold}.} By definition $W^{z^{\ast}}(\cdot,0)=V^{c_{0}}.$ Assuming that $W^{z^{\ast}%
}(\cdot,c_{k})=V^{c_{k}}$ for $k=1,...,i-1$, by Theorem
\ref{Caracterizacion Discreta}, it is enough to prove that $W^{z^{\ast}}%
(\cdot,c_{i})$ is a viscosity solution of (\ref{HJB equation discreta}). Since
by construction $V^{c_{i-1}}-W^{z^{\ast}}(\cdot,c_{i})\leq0$ and $V^{c_{i-1}%
}(x)-W^{z^{\ast}}(x,c_{i})=0$ for $x\leq z^{\ast}(c_{i})$, it remains to be
seen that $\mathcal{L}^{c_{i}}(W^{z^{\ast}})(x,c_{i})\leq0$ for $x\leq
z^{\ast}(c_{i})$. By Remark \ref{Remark Obstaculo}, $W^{z^{\ast}}(\cdot
,c_{i})$ is continuously differentiable and it is piecewise infinitely
differentiable in open intervals in which it solves $\mathcal{L}^{c_{j}%
}(W^{z^{\ast}})(x,c_{i})=0$ for some $j\leq i$. Let us consider first the case
in which $x\neq z^{\ast}(c_{k})$ for $k=1,..,i-1$, so $x$ belongs to one of
these open intervals. Hence,
\[
\mathcal{L}^{c_{i}}(W^{z^{\ast}})(x,c_{i})=\mathcal{L}^{c_{j}}(W^{z^{\ast}%
})(x,c_{i})+(c_{i}-c_{j})(1-\partial_{x}W^{z^{\ast}}(x,c_{i}))=(c_{i}%
-c_{j})(1-\partial_{x}W^{z^{\ast}}(x,c_{i}))\leq0
\]
if and only if $\partial_{x}W^{z^{\ast}}(x,c_{i})\geq1$. Let us prove the
result first for $x=z^{\ast}(c_{i})\neq z^{\ast}(c_{k})$ for $k=1,..,i-1$.
That is, there exists $\delta>0$ and some $j<i$ such that $\mathcal{L}^{c_{j}%
}(W^{z^{\ast}})(x,c_{i})=0$ in $(z^{\ast}(c_{i})-\delta,z^{\ast}(c_{i}))$ and then%

\[
\mathcal{L}^{c_{j}}(W^{z^{\ast}})(z^{\ast}(c_{i})^{-},c_{i})=0\text{,
}\mathcal{L}^{c_{i}}(W^{z^{\ast}})(z^{\ast}(c_{i})^{+},c_{i})=0,
\]
so%
\[%
\begin{array}
[c]{lll}%
0 & = & \mathcal{L}^{c_{i}}(W^{z^{\ast}})(z^{\ast}(c_{i})^{+},c_{i}%
)-\mathcal{L}^{c_{j}}(W^{z^{\ast}})(z^{\ast}(c_{i})^{-},c_{i})\\
& = & \frac{\sigma^{2}}{2}(\partial_{xx}W^{z^{\ast}}(z^{\ast}(c_{i})^{+}%
,c_{i})-\partial_{xx}W^{z^{\ast}}(z^{\ast}(c_{i})^{-},c_{i}))\\
&  & +(c_{i}-c_{j})(1-\partial_{x}W^{z^{\ast}}(z^{\ast}(c_{i}),c_{i})).
\end{array}
\]
By Remark \ref{Remark Obstaculo}, $\partial_{xx}W^{z^{\ast}}(z^{\ast}%
(c_{i})^{+},c_{i})-\partial_{xx}W^{z^{\ast}}(z^{\ast}(c_{i})^{-},c_{i})\geq0$
and $c_{i}-c_{j}>0,$ and we can conclude that $\partial_{x}W^{z^{\ast}}(z^{\ast
}(c_{i}),c_{i}))\geq1$. It remains to prove that $\partial_{x}W^{z^{\ast}%
}(x,c_{i})\geq1$ for $x<z^{\ast}(c_{i}).$

If $i=1$, by definition $W^{z^{\ast}}(x,c_{1})=W^{z^{\ast}}(x,0)=\frac
{\Lambda}{q}\left(  1-e^{\theta_{1}(0)x}\right)  $ for $x\leq z^{\ast}(c_{1}%
)$. By Remark \ref{Remark Concavity}, $W^{z^{\ast}}(\cdot,0)$ is concave and
so $\partial_{x}W^{z^{\ast}}(x,c_{1})=\partial_{x}W^{z^{\ast}}(x,0)\geq
\partial_{x}W^{z^{\ast}}(z^{\ast}(c_{1}),0)\geq1$ for $x\leq z^{\ast}(c_{1})$;
hence we have the result. {We need to prove now that $\partial_{x}W^{z^{\ast}%
}(x,c_{i})=\partial_{x}W^{z^{\ast}}(x,c_{i-1})\geq1$ for $x<z^{\ast}(c_{i})$
and $i>1.$ By the induction hypothesis, we know that $\partial_{x}W^{z^{\ast}%
}(x,c_{i-1})\geq1$ for $x\leq z^{\ast}(c_{i-1})$. In the case that $z^{\ast
}(c_{i})\leq z^{\ast}(c_{i-1}),$ it is straightforward because $x\leq z^{\ast
}(c_{i})\leq z^{\ast}(c_{i-1})$ implies%
\[
\partial_{x}W^{z^{\ast}}(x,c_{i})=\partial_{x}W^{z^{\ast}}(x,c_{i-1})\geq1.
\]
In the case that $z^{\ast}(c_{i})>z^{\ast}(c_{i-1})$, it is enough to prove it
for $x\in(z^{\ast}(c_{i-1}),z^{\ast}(c_{i}))$. Note that $W^{z^{\ast}}%
(x,c_{i})=W^{z^{\ast}}(x,c_{i-1})$ for $x\in(z^{\ast}(c_{i-1}),z^{\ast}%
(c_{i}))$, and $W^{z^{\ast}}(\cdot,c_{i-1})$ is a solution of $\mathcal{L}%
^{c_{i-1}}=0$ in $[z^{\ast}(c_{i-1}),\infty)$ with }$\lim_{x\rightarrow\infty
}{W^{z^{\ast}}(0,c_{i-1})}=\frac{c_{i-1}+\Lambda}{q}${; so, by Remark
\ref{Remark Concavity}, $\partial_{x}W^{z^{\ast}}(x,c_{i})$ is decreasing in
the interval $(z^{\ast}(c_{i-1}),z^{\ast}(c_{i}))$. But $\partial
_{x}W^{z^{\ast}}(z^{\ast}(c_{i}),c_{i})\geq1$, so we have the result. }

Consider now the case $x=z^{\ast}(c_{j})$ with $j\leq i-1$ and $z^{\ast}%
(c_{j})\leq z^{\ast}(c_{i})$. It could be the case that $W^{z^{\ast}}(x,c_{i})$ is not twice
continuously differentiable at $z^{\ast}(c_{j}),$ so we prove that
$\mathcal{L}^{c_{i}}(W^{z^{\ast}})(z^{\ast}(c_{j}),c_{i})\leq0$ in the
viscosity sense. Take a test function $\varphi$ for supersolution at $z^{\ast
}(c_{j}).$ From Definition \ref{Viscosidad Discreta}, $\varphi^{\prime
}(z^{\ast}(c_{j}))=\partial_{x}W^{z^{\ast}}(z^{\ast}(c_{j}),c_{i})$ and%
\[
\varphi^{\prime\prime}(z^{\ast}(c_{j}))\leq\min\{\partial_{xx}W^{z^{\ast}%
}(z^{\ast}(c_{j})^{+},c_{i}),\partial_{xx}W^{z^{\ast}}(z^{\ast}(c_{j}%
)^{-},c_{i})\}.
\]
Then%

\[
\mathcal{L}^{c_{i}}(\varphi)(z^{\ast}(c_{j}),c_{i})\leq\min\{\mathcal{L}%
^{c_{i}}(W^{z^{\ast}})(z^{\ast}(c_{j})^{+},c_{i}),\mathcal{L}^{c_{i}%
}(W^{z^{\ast}})(z^{\ast}(c_{j})^{-},c_{i})\},
\]
and since we already proved that $\mathcal{L}^{c_{i}}(W^{z^{\ast}})(\cdot
,c_{i})\leq0$ in $\left(  z^{\ast}(c_{j})-\delta,z^{\ast}(c_{j})\right)  \cup$
$\left(  z^{\ast}(c_{j}),z^{\ast}(c_{j})+\delta\right)  $ for some $\delta>0,$
we get the result by continuity.\hfill $\blacksquare$\newline

{\color{black} {\bf Data Availability Statement.} No new data were created or analyzed during this study. Data sharing is not applicable to this article.\\	
	
	{\bf Acknowledgements.} The authors would like to thank two anonymous referees and the editor for helpful remarks that improved the presentation of the manuscript. H.A.\ acknowledges support from the Swiss National Science Foundation Project CR00I5-235992. }

\bigskip

\end{document}